\documentclass{amsart}

\usepackage{amsmath, amsfonts, amssymb, amscd}

\newtheorem{theo}{Theorem}[section]
\newtheorem{lem}[theo]{Lemma}
\newtheorem{prop}[theo]{Proposition}

\newtheorem{cor}[theo]{Corollary}
\newtheorem{rem}[theo]{Remark}

\newtheorem{definition}[theo]{Definition}
\newenvironment{pf}{\noindent{\it Proof. }}{$\square$\par\medskip}

\newcommand{\bC}{{\mathbb C}}
\newcommand{\bR}{{\mathbb R}}
\newcommand{\bZ}{{\mathbb Z}}

\newcommand{\cD}{{\mathcal D}}
\newcommand{\cC}{{\mathcal C}}
\newcommand{\cL}{{\mathcal L}}

\newcommand{\cU}{{\mathcal U}}

\newcommand{\cN}{{\mathcal N}}
\newcommand{\cZ}{{\mathcal Z}}

\newcommand{\cH}{{\mathcal H}}

\newcommand{\UU}{{\operatorname{U}}}

\newcommand{\re}{\operatorname{Re}}
\newcommand{\im}{\operatorname{Im}}

\renewcommand{\span}{\operatorname{span}}
\newcommand{\Aut}{\operatorname{Aut}}

\newcommand{\Ker}{\operatorname{Ker}}

\renewcommand{\=}{\overset{\text{def}}{=}}
\newcommand{\codim}{\operatorname{codim}}

\def\sideremark#1{\ifvmode\leavevmode\fi\vadjust{
\vbox to0pt{\hbox to 0pt{\hskip\hsize\hskip1em
\vbox{\hsize3cm\tiny\raggedright\pretolerance10000
\noindent #1\hfill}\hss}\vbox to8pt{\vfil}\vss}}}

\title[Monge-Amp\`ere equations and manifolds of circular type]
{Monge-Amp\`ere equations\\ and moduli spaces of   manifolds of circular type}
\author[G. Patrizio and A.  Spiro]{Giorgio Patrizio
and Andrea Spiro}


\begin{document}

\begin{abstract}   A {\it (bounded) manifold of circular type\/}  is a complex manifold $M$  of dimension $n$ admitting  a (bounded) exhaustive real  function $u$,  defined on $M$ minus a point $x_o$, so that:  a)  it is a smooth solution on $M\setminus \{x_o\}$  to the Monge-Amp\`ere equation $(d d^c u)^n = 0$;  b)  $x_o$ is a singular point for $u$ of logarithmic type and $e^u$ extends smoothly on the blow up of $M$ at $x_o$;  c)  $d d^c (e^u) >0$ at any point of  $M\setminus \{x_o\}$. 
This class  of manifolds naturally includes all smoothly bounded, strictly linearly convex domains and all smoothly bounded, strongly pseudoconvex circular domains of 
 $\bC^n$. \par
The moduli spaces of bounded manifolds of  circular type are studied. In 
particular, for each biholomorphic equivalence class of them it is proved the existence of an essentially unique manifold in normal form.
It is also shown  that the class of normalizing maps for an $n$-dimensional manifold $M$ is a new holomorphic 
invariant with the following property:  it is parameterized by the points of a finite dimensional real manifold of  dimension $n^2$ when $M$ is a (non-convex) circular domain while it is   of dimension $n^2 + 2 n$ when $M$ is  a strictly convex domain. 
 New characterizations of the circular domains 
and of the unit ball  are also obtained.
\end{abstract}

\maketitle

\null \vspace*{-.5in}
\section{Introduction}
\bigskip
In this paper we  analyze
the moduli spaces of a family of complex manifolds, which includes
the smoothly bounded strictly linearly convex domains and the smoothly bounded 
strictly pseudoconvex circular domains in $\bC^n$.  \par
More precisely, we  consider a larger class of manifolds, called
{\it (bounded) manifolds of circular type\/}, which naturally includes the previous two families of domains
and are characterized by the property of  admitting  a (bounded) exhaustive real function $u$, defined on $M$ minus a point $x_o$,  so that:  a)  it is a smooth solution on $M\setminus \{x_o\}$  to the Monge-Amp\`ere equation $(d d^c u)^n = 0$;  b)  $x_o$ is a singular point for $u$ of logarithmic type and $e^u$ extends smoothly on the blow up of $M$ at $x_o$;  c)  $d d^c (e^u) >0$ at any point of  $M\setminus \{x_o\}$.  \par
 In  any biholomorphic equivalence class 
of such  domains, we prove the existence of an essentially unique {\it manifold in normal form\/}, 
consisting of the unit ball $B^n$ together with a non-standard complex structure 
$J$, which satisfies some suitable conditions:  One of them consists on requiring that the  non-standard CR structures induced on 
the spheres  $S^{2n-1}(r) = \{Ê\ |z| = r\ \}$, $ 0 < r < 1$,  have  the same real distribution of  the standard ones  as  underlying distribution of $J$-invariant real subspaces. The other conditions on $J$ imply that
 $J$ is  uniquely determined by only one of such 
  CR structures. This CR structure is  completely determined  by a  sequence $\{\phi_k\}_{k=0}^{\infty}$, 
of $(1,1)$-tensor fields on $\bC P^{n-1} \simeq S^{2n-1}(r_o)/S^1 $,  obtained  by expanding in Fourier series
the tensor field $\phi$  that  gives  the  complex structure of the  CR structure as a deformation  
of the standard one.
As   applications, we use these results to obtain  new characterizations  of the circular domains in $\bC^n$ and of the unit ball.\par
\medskip
The normal forms  considered  in this paper are essentially the same  of the normal forms constructed in \cite{Bl, BD1} by Bland and Duchamp only
for domains that are  small deformations of the unit ball. The major improvement consists in showing  that {\it such normal forms exist  for 
 any bounded manifold of circular type\/}, i.e. for any complex manifold which  admits a solution to the described Monge-Amp\`ere  differential problem.  Such  result has been obtained  by   methods and techniques that are substantially  different from those of  \cite{Bl} and \cite{BD1}. \par
We should also mention that our normal forms  can  be also considered  as  versions  ``without distinguished point'' of  the normal forms constructed 
by  Lempert and  Bland and Duchamp in \cite{Le2, BD},  where equivalence classes of 
 pointed strictly convex domains were studied.\par 
\medskip
In many regards,  the properties of normal forms  of  manifolds of circular type  recall  
 those of the well-known  Chern-Moser  normal forms for   Levi non-degenerate  real  hypersurfaces (\cite{CM}). 
 For instance, if  $D$ is a domain of circular type in a Stein manifold, 
  it turns out that the class $\cN(D)$ of  the diffeomorphisms $f: D \to B^n$,  
which map  $D$  into a   normal form $(B^n, J)$, 
 is naturally parameterized by  a subset of  the automorphism group $\Aut(B^n, J_o)$ 
  of the standard unit ball $(B^n, J_o)$. As for Chern-Moser normal forms, 
 this fact determines   a natural embedding of $\Aut(D)$ as a subgroup of $\Aut(B^n, J_o)$.  \par
 On the other hand,  in contrast with what occurs for
  Chern-Moser  normal forms, the parameter set for the class $\cN(D)$
   is {\it not}  independent of $D$ and it represents  an important biholomorphic invariant for  the domain (we recall that, on the 
   contrary,  the  Chern-Moser normalizing maps are always  parameterized
   by the isotropy $\Aut_{x_o}(B^n, J_o)$  of a fixed boundary point 
$x_o \in \partial B^n$).   \par
For example,   if $D$ is  a smoothly bounded, strictly convex domain in $\bC^n$, then $\cN(D) \simeq \Aut(B^n, J_o)$, 
while if $D$ is a generic  
   (non convex) strongly pseudoconvex circular domain, then $\cN(D) \simeq U_n  \subsetneq \Aut(B_n,J_o)$.  These facts    motivate
  the following question: \par
  \smallskip
  {\it Is it true that a domain of circular type $D$ is biholomorphic to 
  a smoothly bounded, strictly convex domains in $\bC^n$ if and only if    $\cN(D)\simeq \Aut(B^n,J_o)$?\/}
\par
\smallskip
We conjecture that the answer is ``yes'', at least when $n \geq 3$ and some boundary regularity conditions are assumed. However, 
besides such conjecture, there is another reason of interest for the domains of circular type for which  $\cN(D)\simeq \Aut(B^n,J_o)$, namely  the abundance of solutions to the quoted Monge-Amp\`ere problem that there exists on any such domain (there is at least one such solution for any point $x_o \in D$ - see remarks after Thm.  \ref{parameterizationnormalizingmaps}). In addition, 
any such domain is endowed with a biholomorphically invariant complex Finsler metric: In case of   a strictly convex domain of $\bC^n$, this metric is the Kobayashi metric of the domain (e.g. \cite{Fa, Pa, AP, Sp}).  A detailed discussion of this and related questions on such domains will be the object of  a future paper. \par
\medskip
Some  applications of  the theory of normal forms we developed here can be found in the last section. We obtain f.i. the following result (Thm. \ref{lasttheorem}):  {\it a bounded manifold of circular type $M$, with $u: M \setminus \{x_o\} \to [0, r^2)$ satisfying a), b), c),  is biholomorphic to a circular domain in $\bC^n$ if and only if there exists at least two subdomains $M_{< c} = \{ \ u < c\ \}$, $M_{<c'}= \{ \ u < c'\ \}$, $0 < c < c' < r^2$,  which are biholomorphic one to the other  by a map fixing $x_o$.\/}  This theorem 
represents a  generalization to any  manifold of circular type  of results in  \cite{LPW} and \cite{AP},  originally proved  only for smoothly bounded, strictly convex domains in $\bC^n$ or complex Finsler manifolds. Such generalization is obtained by an approach which is   quite different from the ones in \cite{LPW} and \cite{AP} .  By a result in \cite{Pt},  our Thm. \ref{lasttheorem}  has  also an immediate corollary which gives  a new characterization of the unit ball  (Cor. \ref{lastcorollary}).\par
\bigskip
The structure of the paper is as follows: After \S 2,  devoted to preliminaries, in \S 3 we introduce the {\it circular representation\/} of a manifold of circular type, a map which naturally generalizes the standard  circular representation of strictly convex domains (see \cite{Pt, BD}); In \S 4, the normal forms of domains of circular type is defined and the existence of normalizing maps is proved; In \S 5, it is shown that the complex structure of a manifold in normal form is completely determined by the associated ``deformation tensor''  $\phi$;  In the same section, the  Bland and Duchamp invariants $\phi^{(k)}$, determined by Fourier series expansion of $\phi$, are defined; In  \S 6 we give a geometrical interpretation of the Bland and Duchamp invariant $\phi^{(0)}$, we establish the parameterization by elements of $\Aut(B^n)$ of the family of normalizing maps of a manifold of circular type and we prove the mentioned characterizations of circular domains and of the unit ball.\par
 \bigskip 
\section{Preliminaries}
\subsection{Notation, first definitions and some basic properties}\hfil\par
\subsubsection{Complex, CR and contact structures}\hfil\par
\medskip
In all the following, we will be denoted by $J_o$ the standard complex structure of $\bC^n$ 
 and  by  $B^{n}$ and  $S^{2n-1}  = \partial B^n$ the unit  ball and the unit sphere, respectively,  centered at $0 \in \bC^n$. 
 We will also indicate by $\Delta = B^1$ the unit disc in $\bC$.\par
For any  two complex manifolds $(M, J)$ and $(M', J')$,   a map $f: M \to M'$ is  called  {\it  $(J, J')$-holomorphic\/} if $f_* \circ J = J' \circ f_*$. However, 
anytime it will be clear what are the considered complex structures,   we will just write   ``holomorphic"  in place of  ``$(J, J')$-holomorphic"  .\par
\bigskip
 We recall that  a CR structure on a manifold $M$ is  a subbundle $\cD^{1,0} \subset T^\bC M$ of  the complexified 
tangent bundle of $M$ so  that  $\cD^{1,0} \cap \overline{\cD^{1,0}} = \{0\}$ and $[\cD^{1,0}, 
\cD^{1,0}] \subset \cD^{1,0}$. For any given CR structure, we call {\it underlying real distribution\/} the  
  subbundle $\cD \subset TM$ determined by  the subspaces $\cD_x  = \re(\cD^{1,0}_x)\subset T_xM$; We also call
  {\it underlying complex structure\/} the smooth family of complex 
  structures $J_x: \cD_x \to \cD_x$ defined by $J_x(\re(v)) = \re(i \cdot v)$ for any $v \in \cD^{1,0}_x$.
  In the following, a CR structure $\cD^{1,0}$ will be often indicated by  the associated pair $(\cD, J)$, 
of underlying  distribution  and   complex structure.   Notice that any CR distribution $\cD^{1,0}$ can be completely recovered from its associated pair  $(\cD, J)$: 
In fact,  at any $x \in M$, the subspace $\cD^{1,0}_x \subset T^\bC M$ coincides with the $J_x$-eigenspace  in $\cD^\bC_x$ of eigenvalue $+i$.\par
\bigskip
A {\it CR-equivalence\/} between two CR manifolds $(M, \cD^{1,0})$ and  $(M', \cD'{}^{1,0})$ is a diffeomorphism 
$\varphi: M \to M'$ such that $\phi_*(\cD^{1,0}) = \cD'{}^{1,0}$. If we 
consider the pairs  $(\cD, J)$ and  $(\cD', J')$  associated with $\cD^{1,0}$ and $\cD'{}^{1,0}$, respectively, it follows immediately from 
definition that  a diffeomorphism $\varphi: M \to M'$
is a   CR-equivalence if and if $\phi_*(\cD) = \cD'$ and $\phi_*(J) = J'$
\par
\bigskip
Let $(\cD, J)$ be a CR structure of hypersurface type (i.e. with $\codim \cD =  1$) on 
a manifold $M$  with  $\dim M = 2n-1$. Such CR structure is said to be  {\it Levi non-degenerate\/}
if  the underlying real distribution  $\cD$ is a {\it contact distribution\/}. This corresponds to say that    any  1-form
$\theta$, satisfying   $\Ker \theta|_x = \cD_x$ at any $x$,  is a contact form, 
i.e. so that $\theta \wedge (d\theta)^n  \equiv 0$ or, equivalently, so that 
$d\theta_x|_{\cD_x \times \cD_x}$ is  a non-degenerate 2-form on $\cD_x$.  We would like to point out that this definition
 is completely equivalent to the classical definition of  ``Levi non-degeneracy''.\par
\par
\bigskip
For any real hypersurface $S \subset \bC^n$,  the CR structure of hypersurface type  induced on $S$ from 
the standard complex structure $J_o$ will be denoted by $(\cD_o, J_o)$ and it will be called  {\it standard CR structure of $S$\/}. 
When $S$ is the  smooth boundary of a strongly pseudoconvex domain, the 
distribution $\cD_o$  is a a contact distribution and  will be called {\it the standard contact distribution of $S$\/}.
\par
\bigskip
\subsubsection{Circular domains and Minkowski functions}\hfill\par
\medskip
For any $\zeta\in \bC$,  let us  denote by $\zeta \cdot (\ \ ): \bC^{n} \to \bC^{n}$
the holomorphic transformation 
$$\zeta\cdot (z^1, \dots, z^{n+1}) \=  (\zeta \cdot z^1, \dots, \zeta \cdot z^{n+1}) \ . \eqno(2.1)$$
The  circle $S^1 = \partial \Delta$ will be often
tacitly identified  with the  compact 1-parameter group of holomorphic  transformations $S^1 = \{\ e^{i t}\cdot(\ )\ \}$.
We recall that a domain $D\subset \bC^{n+1}$ is called {\it circular\/} if
it is invariant under any transformation  of $S^1$, while it is 
called {\it circular and complete\/} if it is   invariant under all  
transformations $\zeta \cdot (\ \ )$ with 
 $\zeta \in \bar \Delta$. \par
\medskip
For any  complete circular domain  $D \subset \bC^n$,  the associated 
 {\it Minkowski function \/} is the map $\mu_D: \bC^{n+1} \to \bR_{\geq 0}$ defined by 
$$\mu_D(z) = \frac 1 t_z\ ,\ \qquad  \text{where}\ t_z = \sup\{\ s\in \bR\ :\ s\cdot z\in D\ \}\ .\eqno(2.2)$$
Notice that,   $\mu_D(\zeta \cdot z) =  |\zeta| \mu_D( z)$ and 
that $D = \{\ z\in  \bC^{n+1}\ :\ \mu_D(z) < 1\ \}$.
In particular,  if $D$ is  smoothly bounded,  the function $\tau = \mu_D^2 - 1$ is a  defining 
functions for $D$ which is smooth on $\bC^{n+1} \setminus\{0\}$.  In this case,  
$D$ is  strictly linearly convex if and only if    $\tau = \mu_D^2 - 1$  has strictly positive Hessian  at any  $x \neq 0$.\par
\smallskip
We conclude with the following definition. If $D\subset \bC^n$ is a complete circular domain with Minkowski function $\mu_D$, for any $v \in T_0 \bC^n = \bC^n$ we call {\it standard radial disc of $D$ tangent to $v$\/} the holomorphic map
 $$f^{(\mu)}_v: \bar \Delta \to D \ ,\qquad f^{(\mu)}_v(\zeta) = \zeta \cdot \frac{v}{\mu(v)}\ .$$
 For instance, in case $D = B^n$  the ``standard radial discs'' are  the holomorphic discs  of the form  $f_v(\zeta) = \zeta \cdot \frac{v}{|v|}$ for some $v \in \bC^n$. We remark that 
 any standard radial disc of a smoothly bounded complete circular domain $D$  is indeed a {\it stationary disc for $S = \partial D$\/}, according to the well-known definition of  Lempert (\cite{Le}; for a definition of stationary discs of an hypersurface see also \cite{Tu}); for a proof of such property we refer  e.g.  \cite{Pa1}, Lemma 3.34. \par
 \bigskip
 \subsubsection{Blow ups of  $\bC^n$ and lifts of standard radial discs: some elementary facts}\hfill\par
 \medskip
In all this paper,  we will denote by $\widetilde \bC^n$ the blow up at the origin of $\bC^n$, 
by $\widetilde \pi: \widetilde \bC^n \to \bC^n$ the standard projection and, for any domain $D \subset \bC^n$ containing 
the origin, 
we will indicate  by $\widetilde D$  the  blow up of $D$ at $0$,  considered as a subset of $\widetilde \bC^n$. In other 
words,   $\widetilde D\= \widetilde \pi^{-1}(D)\subset \widetilde \bC^n$. \par
We 
recall that $\widetilde \bC^n$ coincides with the tautological 
line bundle $\pi: \widetilde \bC^n = E \longrightarrow \bC P^{n-1}$, i.e. with the bundle $E$ over $\bC P^n$
given by the  pairs 
$([v], z) \in \bC P^{n-1} \times \bC^n$ such that  $z \in [v]$, and that  the exceptional divisor of $\widetilde \bC^n = E$ 
coincides with the image of  the zero section of $\pi: E \to \bC P^n$.\par
\smallskip
Given  a system of 
affine coordinates $v =  (v^1, \dots, v^{n-1}): \cU \subset \bC P^{n-1} \to \bC^{n-1}$ for some open subset of  $\bC P^{n-1}$
we call {\it associated system of coordinates\/} the system of 
 coordinates $(v, \zeta) = (v^1, \dots, v^{n-1}, \zeta):  \widetilde \cU = \pi^{-1}(\cU)\subset E \to \bC^n$,  which maps  any point $([v], z) \in \widetilde \cU$
 into the $n$-tuple, whose first $n-1$ components are 
 the affine coordinates of $[v]$, while   the last is 
  the unique complex number so that 
$$([v], z) = \left([v^1, \dots, \underset{\text{i-th\ place}} 1,\dots v^{n-1}],   \zeta\cdot v^1, \dots, \underset{\text{i-th\ place}} \zeta,\dots ,\zeta \cdot v^{n-1})\right)\ .\eqno(2.5)$$
\par
\bigskip
We conclude this short subsection noticing that  any standard radial disc $f_v: \Delta \to B^n$
admits a unique lifted map 
$$\widetilde f_v: \Delta \to \overline{\widetilde B^n}\ ,\eqno(2.6)$$
which projects down to $f_v$. It can be checked that any map $\widetilde f_v$
 is a  stationary disc for the sphere $S^{2n-1}$, considered this time as a real hypersurface of 
$\widetilde \bC^n$. Moreover, the image of $\widetilde f_v$  coincides with one of the fibers of $\pi_{B^n}: \widetilde B^n \to \bC P^{n-1}$.
In the following, we will shortly refer to the discs $\widetilde f_v$ as the {\it standard radial
discs of $\widetilde B^n$}: Notice that  the images of such standard radial discs 
determines a holomorphic foliation for $\widetilde B^n$ (they are the fibers of the bundle structure 
of $\widetilde B^n$ over $\bC P^{n-1})$. \par
 Analogous conclusions hold for the lifts in $\widetilde \bC^n$ of the the radial discs $f^{(\mu)}_v$ 
 of a  complete circular domain $D$.\par
\bigskip
\subsection{Manifolds of circular type, indicatrices and Monge-Amp\`ere foliations}\hfil\par
\medskip
In this section, we give the definition of  
``manifolds of circular type'', a notion introduced by the first author in  \cite{Pt1}, and we  define  a 
few other related concepts that will reveal to be essential in the study of the moduli spaces of such 
manifolds. In particular, we are going to show in the next section that 
{\it any manifold of circular type admits a ``circular representation''\/} (see \S 3 for definition) 
extending a  property usually stated only for strictly convex domains or circular domains in $\bC^n$.  \par
\subsubsection{Manifolds and domains of circular type}\hfil\par
\medskip
In the next phrases, for any function $\tau: M \to \bR$, 
we will denote  $M_{\tau = 0} = \tau^{-1}(0)$ and
$M_{\tau \neq 0} = M \setminus M_{\tau = 0}$.
\par
\begin{definition} \label{circulartype} {\rm Let $M$ be a non-compact complex manifold of complex dimension 
$n$. We say that $M$ is  a {\it  manifold of circular type\/} if  it admits an exhaustion function 
$\tau: M \longrightarrow [0, r^2)$, for some $r^2 \in (0, \infty]$, such that
\begin{itemize}
\item[a)] $\tau \in \cC^0(M) \cap \cC^\infty(M_{\tau \neq 0})$;
\item[b)] on $M_{\tau \neq 0}$, $\tau$ is so that $d d^c \tau > 0$ and $d d^c \log \tau \geq 0$; 
\item[c)] on $M_{\tau \neq 0}$, $(d d^c \log \tau)^n \equiv 0$;
\item[d)] there exists a point $x_o \in M_{\tau = 0}$ so that, if $\pi: \widetilde M \to M$ is the blow up
of $M$ at $x_o$, then $\tau \circ \pi: \widetilde M \to \bR$ is smooth and there exist  two positive constant 
$C_1, C_2$
so that $C_1 \|x - x_o\|^2 \leq \tau(x) \leq C_2\|x - x_o\|^2 $ at all points of a neighborhood of $x_o$
(here ``$\|\cdot \| $" is the Euclidean norm in some coordinate neighborhood centered at $x_o$).
\end{itemize}
It is easy to see that $\{ \tau = 0\} = \{x_o\}$ (\cite{Pt1}). Such point $x_o$ is called  {\it center of $M$ w.r.t. $\tau$\/} and   $\tau: M \longrightarrow [0, r^2)$ 
is called  {\it parabolic exhaustive function\/}. If $r < \infty$ (i.e. $\tau$ is bounded), we call  $M$  {\it bounded\/}.
\par
A smoothly bounded, relatively compact  domain $D$ in a complex manifold $(M,J)$ is called {\it domain of circular type} if $(D, J)$
 is a bounded manifold of circular type admitting at least one  parabolic exhaustive function  $\tau$, which admits a smooth extension 
up to the boundary and satisfies b) also at the  points of $\partial D$.  For domains of circular type, we shall consider only exhaustions that are smooth  and satisfying b) up the boundary. For such domains, the expression  {\it parabolic exhaustion\/} will always mean such a function}.
\end{definition}

\begin{lem} \label{uniqueness} \hfill
 \begin{itemize}
\item[a)]   Any manifold of circular type  is Stein;
\item[b)]  if $(M, J, \tau)$ is a domain of circular type and  $\tau, \tau'$  are two parabolic exhaustive functions of $M$ (smooth up to the boundary)  having the same center $x_o$, then  $\tau' = k \tau$ for some suitable positive constant $k$.
\end{itemize}
\end{lem}
\begin{pf}  Claim (a) is immediate (see f. i. \cite{Na}).
Claim (b) is  a consequence of Thm. 4.3 in \cite{De}. In fact,   $M$ 
is  a hyperconvex domain according to  Def. 2.1 in \cite{De} and hence, for the cited theorem, there exists a unique 
plurisubharmonic function $u$ so that $u|_{\partial M} = 0$, $(d d^c u)^n \equiv 0$ and $u(z) \sim \log|z - x_o|$ 
for $z \to x_o$. This implies that $\log \tau = u + c$ and $\log \tau' = u + c'$ where $c$ and $c'$ are the constant values  $c = \log \tau|_{\partial M}$, 
$c' = \log \tau|_{\partial M}$. From this the claim follows immediately.
\end{pf}
\medskip
Any complete circular domain of $\bC^n$ is a domain of circular type.  
In fact, for any such domain $D$,  if we denote by $\mu_D: \bC^n \to \bR_{\geq 0}$  the Minkowski functional of $D$, 
then  the map $\tau = \mu^2_D$ is a  parabolic exhaustive function   with center $x_o = 0$. Notice  that,  {\it generically, 
the origin $0$ 
 is the only center  of  a complete circular domain\/}.\par
\medskip
Other examples of domains of circular type are given by all bounded,
strictly convex domain  $D \subset \bC^n$ with smooth boundary. In fact, using the  properties 
of the Kobayashi pseudo-distance of such domains (see \cite{Le, Le1, Pt1}), we have that 
 for any point $x_o \in D$, the function
$$\tau_{x_o}: D \longrightarrow \bR_+\ ,\qquad \tau_{x_o}(x) =  \left(\frac{e^{ \kappa_D(x,x_o) - 1}}{ 
e^{ \kappa_D(x,x_o) + 1}}\right)^2\ ,$$
($\kappa_D = $  Kobayashi pseudo-distance of $D$)  is 
a parabolic exhaustive function for $D$ which extends smoothly up to the boundary. In particular, {\it any point $x_o$ of a smoothly bounded, strictly convex domain $D$ in $\bC^n$  
  is a center for  $D$ w.r.t. some parabolic exhaustive function\/}. \par
 \bigskip 
 \subsubsection{Monge-Amp\`ere foliations and normal distributions}\hfill\par
\medskip
Let $(M, J, \tau)$ be a bounded manifold of circular type and $x_o \in M$ be the 
center for $M$ w.r.t.  $\tau$.  By multiplication of $\tau$
by positive constant, we may  assume that $\tau: M \to [0, 1)$.  We will also denote by $\widetilde M$ the blow up 
of $M$ at $x_o$.\par
\smallskip
First of all, we want recall a basic  identity that follows from property (c) of Definition \ref{circulartype}
 (see e.g. \cite{St}).
Using just the definitions, one finds immediately that  
$$ \tau^2 d d^c \log \tau =\tau d d^c \tau - d \tau \wedge d^c \tau \ , \eqno(2.7)$$
$$ \tau^{k+1} (d d^c \log \tau)^k =\tau (d d^c \tau)^k - k d \tau \wedge d^c \tau \wedge (d d^c \tau)^{k-1}\ .\eqno(2.8)$$
So,  condition (c) of Definition \ref{circulartype} is equivalent to the equality
$$ \tau (d d^c \tau)^n = n d \tau \wedge d^c \tau \wedge (d d^c \tau)^{n-1}\eqno(2.9)$$
that has to be satisfied at all points of $M \setminus \{x_o\}$. 
\par
\bigskip
Consider now the vector field $Z$ on $M \setminus \{x_o\}$ defined by the 
condition
$$d d^c \tau(Z, J X) =  X(\tau) \ \qquad \text{for any}\ \ X \in T\left( M \setminus \{x_o\}\right)\ .\eqno(2.10)$$
Let also denote by $\cZ$ and $\cH$ the $J$-invariant distributions on $M \setminus \{x_o\}$ defined
as the family of spaces \par
\smallskip
$$\cZ_x \= \span_\bR\{\ Z_x\ ,\ J Z_x\ \}\ ,\eqno(2.11)$$
$$ \cH_x = (\cZ_x)^\perp \= \{\ X \in T_xM\ :\ dd^c\tau(Z, X) = dd^c\tau(J Z, X) = 0 \ \}\ .\eqno(2.12)$$
\par
\smallskip
\noindent 
and 
notice that  $T_x M = \cZ_x \oplus \cH_x$.  Using (2.10) and  (2.9), it follows  that
$$\left.d d^c \tau(Z, JZ)\right|_x = \left.d\tau(Z)\right|_x = \tau_x\ .\eqno(2.13)$$
From  (2.13) and (2.7), it  can be checked  that $\cZ$
coincides with the family of subspaces 
 $$\cZ_x = \ker \left.d d^c \log \tau\right|_x\eqno(2.14)$$ 
 and  hence  {\it is an 
integrable distribution\/}. 
On the other hand, from (2.12) and  (2.10), it turns out that  any space $\cH_x \subset \cH$
 is tangent to  the real hypersurface $S_x = \{\ y\ : \tau(y) = \tau(x)\ \}$. By the fact that $\cH_x$  is $J$-invariant and of real dimension $2(n-1)$, it follows that
 {\it for any $S_c \= \{\ \tau = c\ \}$, the distribution $\cH|_{S_c}$ is  the real distribution underlying the CR structure of  $S_c$\/}.\par
 The distributions $\cZ$ and $\cH$ are defined just on  $M \setminus \{x_o\}$ and they  do not  extend to a regular distribution 
 over the whole $M$. But they do extend on $\widetilde M$: 
 \medskip
 \begin{lem} \label{smoothextendibility}The distributions $\cZ$ and $\cH$ extend  in a unique way as  smooth $J$-invariant distributions  on 
the whole blow up $\widetilde{M}$. The extension of $\cZ$ is integrable everywhere, 
while the extension of $\cH$ is   integrable only when restricted on $\pi^{-1}(x_o)$. More precisely, after identification of a neighborhood of $\pi^{-1}(x_o) \subset \widetilde M$
with an open neighborhood of $\bC P^{n-1}$ in $E = \widetilde \bC^n$, 
then $\cZ|_{\bC P^{n-1}}$ 
coincides with the restriction to $\bC P^{n-1}$ of the vertical distribution of $\pi: E\to \bC P^n$, while $\cH|_{ \bC P^{n-1}} = T  \bC P^{n-1}$.
 \end{lem}
 \begin{pf} For what concerns the extendibility of $\cZ$, it suffices to re-write   formulae (2.7) and (2.10) of \cite{Pt1}
 using a system of coordinates $(v^1,\dots, v^{n-1},  \zeta)$ as described in (2.5).  From those formulae
 it is immediate to realize that  the vector fields $Z$ and $J Z$ extend  smoothly at   $\bC P^{n-1} = \pi^{-1}(x_o) \subset \widetilde M$ 
  setting $Z_{([v],0)} = \re\left(\left.\partial/\partial \zeta\right|_{([v],0)}\right)$, $J Z_{([v],0)}  = \im\left(\left.\partial/\partial \zeta\right|_{([v],0)}\right)$, 
  and that the extension of $\cZ$  is  integrable. \par
About the extendibility of  $\cH$, we recall that (see e.g. \cite{Wo})
$$\cL_Z d d^c \tau = d\left(\imath_Z d d^c \tau\right) \overset{(2.19)}= $$
$$ = d\left(\frac{1}{\tau} \imath_Z d \tau \wedge d^c  \tau + \tau\cdot  \imath_Z d d^c \log \tau\right) \overset{(2.12), (2.8), (2.11)}=  d d^c \tau\ .\eqno(2.15)$$
From (2.15) and (2.12), it follows that   $\cH$ is preserved by  the flow  $\Phi^Z_t$ of $Z$ on $M \setminus \{x_o\}$\footnote{The same  arguments
imply that $\cH$ is preserved also by the flow of $J Z$.} 
 and hence 
it is natural to extend $\cH$  by setting,  for any $y \in \bC P^{n-1} = \pi^{-1}(x_o) \subset \widetilde M$, 
$$\cH_{y}  = \Phi^Z_{-t}{}_*(\cH_{\Phi^Z_t(y)})\ \qquad \text{for some sufficiently small}\ \  t > 0\ .\eqno(2.16)$$
It is quite direct to check that (2.16)  does not depend on $t$ and that in this way $\cH$  is extended 
smoothly.  Moreover, by Lemma 2.1 in \cite{Pt1},  there exists a smooth positive map $h$ such  that,   in a system of coordinates (2.5),   $\tau$ is of the form
$$\tau(v, \zeta) = |\zeta|^2 h^2(v) + o(|\zeta|^3)\eqno(2.17)$$ 
(see (2.18) below). So,  for any $y  = ([v], 0) \in \bC P^{n-1}$, $\cH_y = T_y \bC P^{n-1}$.
\end{pf}
The integrable distribution $\cZ$ and its associated foliation, usually called {\it Monge-Amp\`ere foliation\/}, were 
 considered for the first time by Bedford and Kalka in \cite{BK}. We will therefore call 
 the complementary distribution $\cH$ the {\it normal  distribution of the  Monge-Amp\`ere foliation\/}.\par
\bigskip
\begin{lem} \label{seconduniqueness} Let   $(M, J, \tau)$ be a manifold of circular type and  $\tau, \tau'$ two   parabolic exhaustive functions of $M$ with the same associated  center $x_o$. If there exists $c, c'> 0$ such that $ \emptyset \neq \{\ \tau = c\ \} =
\{\ \tau' = c'\ \}$, then   $\tau' = k  \tau$ for some positive constant $k$.
\end{lem}
\begin{pf} By Lemma \ref{uniqueness}, if we denote by $D = \{\ \tau < c\ \} = \{\ \tau' < c'\ \}$, then 
$\tau|_D = k \tau'|_D$ for some $k>0$. Replacing $\tau'$ by $k \cdot \tau'$, it remains to show that 
if  $\tau|_D = \tau'|_D$  then  $\tau = \tau'$.  For this,  let us denote by  $\widetilde \tau = \tau \circ \pi$ and $\widetilde \tau' = \tau' \circ \pi$ the lifts  of $\tau, \tau'$ at the blow up $\pi: \widetilde M\to M$. Let also $\widetilde D = \pi^{-1}(D)$.  Consider a leaf $S$ of the Monge-Amp\`ere foliation of $\widetilde M$, determined by $\tau$,  and observe that $S \cap \widetilde D$ coincides with an the open subset $S'\cap \widetilde D$ of some leaf $S'$ of  the Monge-Amp\`ere foliation determined by $\tau'$. 
Being $S$ and $S'$ analytic submanifolds of $\widetilde M$,  $S = S'$ and the Monge-Amp\`ere foliations of $\tau$, $\tau'$ coincide.  Now,  let  $Z$ and $Z'$ be the vector fields defined via (2.10) by $\tau$ and $\tau'$, respectively,  and consider the  complex vector fields $Z^{1,0}= \frac{1}{2}(Z - i J Z)$,   $Z'{}^{1,0}= \frac{1}{2}(Z' - i J Z')$. The restrictions $Z^{1,0}|_S$ and $Z'{}^{1,0}|_S$ on a leaf $S$ of the Monge-Amp\`ere foliation are holomorphic vector fields for $S$. Since  they coincide  on $S \cap \widetilde D$, we get  that $Z^{1,0}|_S \equiv Z'{}^{1,0}|_S$. By arbitrariness of $S$,  $Z = Z'$ and hence $\tau = \tau'$ by (2.13).
\end{pf}  
\medskip
We conclude this section with the following technical lemma, which will be needed in the following 
and shows that $\cH$ is uniquely determined by the center $x_o$ and  the  vector field  (2.10). \par
Let $(M, J, \tau)$ be a bounded manifold of circular type 
and $\pi: \widetilde M \to M$ the usual blow up at the center $x_o$. Let also $(J', \tau')$ 
 another pair   so that also   $(M, J', \tau')$ is a domain of circular type with the same center 
$x_o$ and let $\pi: \widetilde M' \to M$ the new blow up at $x_o$. Then, we have the following:\par
\begin{lem} \label{normaldistribution}Ê
If $(J, \tau)$ and $(J', \tau')$  determine the same vector field $Z$  via (2.10), then $\widetilde M = \widetilde M'$ (as real manifolds) , 
$\tau = \tau'$ and   $\cH = \cH'$, i.e. the normal distributions of  the Monge-Amp\`ere 
foliations of $(M, J, \tau)$ and $(M, J', \tau')$  are the same.\par
In particular, the CR structures induced by  $J$ and $J'$ on any  hypersurface $S_c = \{\tau = \tau' =  c\}$ have the same underlying real distributions.
\end{lem}
\begin{pf} Consider the  identity map between $\widetilde M \setminus \pi^{-1}(x_o) = M \setminus \{x_o\} = \widetilde M' \setminus \pi'{}^{-1}(x_o)$ and extend it continuously along each integral curve of $Z$. 
From  Lemma \ref{smoothextendibility}, 
such map is  unique and smooth relatively to the real manifold structures of  $\widetilde M$ and $\widetilde M'$ and implies  that we may consider $\widetilde M= \widetilde M'$ as real
differentiable manifolds. 
To see that $\tau = \tau'$, it suffices to observe that, by (2.17) and  (2.14), $\tau|_{\pi^{-1}(x_o)} = \tau'|_{\pi^{-1}(x_o)} = 0$ and, for any   $x$ of the form $x = \Phi^Z_t(y)$,  $t > 0$, for some  $y \in \pi^{-1}(x_o)$,  we have that
$ \tau(x) = e^t  = \tau'(x)$. Finally, since by Lemma \ref{smoothextendibility}
$\cH_y = \cH'_y = T_y \bC P^{n-1}$ for any 
point $y \in  \bC P^{n-1} = \pi^{-1}(x_o)$, by the invariance of $\cH$ and $\cH'$ under the 
flow of $Z$, the same argument of before implies that $\cH_x = \cH'_x$ at any $x \in M$.\end{pf}
\bigskip
\subsubsection{Indicatrices  of a  manifold of  circular type}\hfil\par
\medskip
Let $(M, J, \tau)$ be a bounded manifold  of circular type and $x_o \in M$ the 
center for $M$ w.r.t.  $\tau$. As before, we assume that $\tau: M \to [0, 1)$. 
\par
\smallskip
Consider a system of complex coordinates $(z^1, \dots, z^n)$ on a neighborhood of $x_o$
with $z^j(x_o) = 0$,  $j = 1, \dots, n$.  By Lemma 2.1 in \cite{Pt1},  there exists a smooth map $h: S^{2n-1} \to \bR_{>0}$ such 
that 
$$\sqrt{\tau(z)} = |z| h\left(\frac{z}{|z|}\right) + o(|z|^2)\ . \eqno(2.18)$$
We set 
 $$\kappa:  T_{x_o} M \simeq \bC^n \to \bR_{\geq 0}\ ,\qquad \kappa(v) =
 \left\{\begin{matrix} |v| \cdot h\left(\frac{v}{|v|}\right)  &\text{if}\ v \neq 0\\
 0 & \text{if} \ v = 0
 \end{matrix}\right.\ .\eqno(2.19)$$
  It can be immediately 
  checked that the value $\kappa(v)$ coincides with 
 $$\kappa(v) = \lim_{t_o\to 0} \left.\frac{d }{dt} \sqrt{\tau(\gamma_t)}\right|_{t = t_o}\ ,$$
 for any smooth curve $\gamma_t$ with $\gamma_0 = x_o$ and $\dot \gamma_0 = v$. In particular,
 it  does not  depend on the choice of the coordinate 
 system. Furthermore 
 we have that 
 $$\kappa (\lambda v) = |\lambda| \kappa(v)\qquad \text{for any} \  \lambda \in \bC$$
 and that $\kappa$ is a smooth function on  $T_{x_o} M \setminus \{0\}$. \par
\begin{definition}{\rm  The {\it indicatrix of $M$ at the center $x_o$ determined by $\tau$\/} is the  smoothly bounded, complete circular domain $I_{x_o} \subset T_{x_o} M$ defined by 
$$I_{x_o} = \{\ v \in T_{x_o} M \ :\  \kappa(v) < 1\ \}\ .\eqno(2.20)$$
}
\end{definition}
\begin{rem} {\rm Notice  that $I_{x_o}$ and $\kappa^2$ coincide with the domain $G(r)$ and the function $\sigma$ defined in (4.3)  and  (4.1) of \cite{Pt1}, respectively. In particular, 
 by the proof of Prop. 4.1 of that paper,  $d d^c \kappa^2 >0$ at all points of $\overline I_{x_o}\setminus \{0\}$ and hence 
 $I_{x_o}$ is  strongly pseudoconvex. Observe also that 
 if  $M$ is a strictly convex domain in $\bC^n$ and $\tau$ is as in (1.1), then  $\kappa$ coincides 
 with the infinitesimal Kobayashi metric of $M$ at $x_o$ and  $I_{x_o}$ is  the 
 Kobayashi indicatrix of $M$ at $x_o$.}
 \end{rem}
\bigskip
\section{Circular representations}
\bigskip
In this  section,   $(M, J, \tau)$ is a  bounded manifold of circular type, $\tau$ takes values in $ [0, 1)$,   $x_o \in M$  is the center associated  with $\tau$ and $I  = I_{x_o}\subset T_{x_o}M$ the corresponding indicatrix at $x_o$. 
Keep in mind that, under an identification $T_{x_o} M \simeq \bC^n$, $I$ can be considered as a circular 
domain, whose Minkowski function is the map $\kappa$ defined in (2.19). We will also denote by $\widetilde M$ the blow up 
of $M$   at $x_o$  and we systematically insert the tilde ``\ $\widetilde{\phantom{a}}$\ "  on top of symbols of manifolds,  domains or maps,  whenever 
 we want to indicate a blow up    or a lift of a map at  such blow up. In particular, $\widetilde \tau = \tau \circ \pi: \widetilde M \to [0,1)$, 
  $\widetilde I$ and $\widetilde{T_{x_o} M}$ are  the blow ups  of $I$ and $T_{x_o} M \simeq \bC^n$ at $0$ and, for any $v \in \partial I$, 
    we set   $\widetilde f_{v}: \Delta \to \widetilde I$ to be the 
the radial   disc  of $\partial \widetilde I $,  obtained by lifting the radial disc $f^{(\kappa)}_v$  and defined in (2.6) and following. 
\par
\medskip
 The purpose of this section is to  recall a few well known properties, usually stated for strictly convex domains 
 (see e.g. \cite{BD, BDK, Le, Le2, Pa1, Pt}) and 
which actually hold for all  domains of circular type. 
We collect  such properties in the next proposition, whose  proof 
is   a direct application of the results in  \cite{Pt1}.  Since the notation and terminology used \cite{Pt1} is 
substantially different from  the one of this paper,  we give also an outline of the proof.
\smallskip
\begin{prop} \label{circularrepresentation} Let $(M, J, \tau)$ be a bounded domain of circular type with 
 $\tau: M \to [0,1)$, 
and let $I \subset T_{x_o} M$ be the indicatrix 
determined by  $\tau$.  Then  
  there exists a unique diffeomorphism $\Psi: \widetilde{I} \to \widetilde M$ 
with the following properties: 
\begin{itemize}
\item[i)]  $\Psi|_{\pi^{-1}(0)} = Id_{\pi^{-1}(0)}$,  provided that we naturally identify
the exceptional divisor  $\pi^{-1}(0)$  of  $\widetilde{ T_{x_o} M}$ with  
 the exceptional divisor  $\pi^{-1}(x_o)$ of $ \widetilde M$; 
\item[ii)]ÊFor any $t \in ]0,1[$,  consider the map
$$\Phi^{(t)}: \partial \widetilde I \to \widetilde M\ ,\qquad \Phi^{(t)}([v], z) = \Psi([v], t z)\ ;\eqno(3.1)$$
 Then  $ \Phi^{(t)}|_{ \partial \widetilde I}$ is a diffeomorphism between 
 $ \partial \widetilde I$ and the level hypersurface $S^{(t)} =  \{\  \tau=  t^2\ \} \subset \widetilde M$; 
 \item[iii)]Êfor any   $t \in ]0,1[$, $ \Phi^{(t)}|_{ \partial \widetilde I}$ 
 maps  the real distribution of  the CR structure of $\partial \widetilde I$ onto the real distribution 
of  the CR structure of   $S^{(t)}$; 
\item[iv)] for any    $([v], z) \in \partial \widetilde I$ and $t \in ]0, 1[$,   the map 
$$\widetilde f^{(t)}_{([v], z)}  :\Delta \to \widetilde M\ ,\qquad \widetilde f^{(t)}_{([v], z)}( \zeta) = \Psi([v], t \zeta z)\eqno(3.2)$$ 
is holomorphic,  injective, with $\widetilde f^{(t)}_{([v], z)}(\partial \Delta) \subset S^{(t)}$ and with 
$\widetilde f^{(t)}_{([v], z)}(\Delta)$equal to an integral leaf of the Monge-Amp\`ere 
foliation of $\widetilde M_{< t^2} = \{\ \widetilde \tau < t^2\ \}$. 
\end{itemize}
In case $M$ is a domain of circular type and $\tau$ is smooth up to the boundary, then $\Psi$ extends 
smoothly to the boundary  with 
$\Psi(\partial \widetilde I) \subset  \partial \widetilde M$ and (iii), (iv) are valid also for $t = 1$.
\end{prop}
\begin{pf}  Pick a  complex basis for $T_{x_o} M$ and use it 
to   identify  $T_{x_o} M \simeq \bC^n$ in the rest of the proof. Then,  consider the map
$$\Psi: I \setminus \{0\} \longrightarrow M\setminus\{x_o\} \ ,\qquad \Psi(v) = F\left( \frac{v}{|v|}, \kappa(v)\right)\eqno(3.3)$$
where $F: S^{2n-1}  \times \Delta \to M$ is the  smooth map  of Thm. 3.4 of  \cite{Pt1}, with   $S^{2n-1}$   unit sphere in $\bC^n =  T_{x_o} M$ and $\Delta$ unit disc in $\bC$. 
In \cite{Pt1} it is proved that:
\begin{itemize}
\item[a)] for any $v \in S^{2n-1}$,  the map  
$f^{(v)} : \Delta \to M$, defined by   $f^{(v)}(\zeta)  
= F\left( v ,  \zeta \right)$ is proper, holomorphic and injective,    with $f^{(v)}(0) = x_o$ and such that $f^{(v)}(\Delta \setminus \{0\})$ is 
a leaf of  the Monge-Amp\`ere foliation;
\item[b)] if  $(\rho, \theta)$ are the standard polar coordinates of $\bC$, the map in (a) is so that 
$$f^{(v)}_*\left(\left.\frac{\partial}{\partial \rho}\right|_\zeta\right) = 2 \left.\sqrt{\tau}\cdot Z\right|_{F(v, \zeta)}\ ,\qquad 
\text{for any}\ \zeta \in \Delta \setminus \{0\}\eqno(3.4)$$
$$f^{(v)}_*\left(\left.\frac{\partial}{\partial \rho}\right|_0\right) = \frac{v}{\kappa(v)} \in T_{x_o} M\ \ ;\eqno(3.5)$$
\item[c)] $F(e^{i \theta} \cdot y, \zeta) = F(y, e^{i \theta} \zeta)$ for any $\theta \in \bR$;
\item[d)] $\tau(F(v, \zeta)) = |\zeta|^2$ for any $v \in S^{2n-1}$.
\end{itemize} 
Now, identify  the exceptional divisor  of  $\widetilde{ T_{x_o} M}$ with  
 the exceptional divisor   of $ \widetilde M$  and extend $\Psi$ to a map 
  $\Psi: \widetilde I \to \widetilde M$  by  setting $\Psi|_{\bC P^{n-1}} = Id_{\bC P^{n-1}}$. 
 Using definitions, Thm. 3.4 and Lemma 2.5 of \cite{Pt1}, 
it is possible to check that such  extended  map
 is a diffeomorphism.\par
Claim (ii) follows directly   from d) and  the definition of  $\Psi$. To check (iii), notice that 
 $\Phi^{(t)} = \Psi \circ \delta^{(t)}$ where $\delta^{(t)}$ is the contraction
$$\delta^{(t)}: \partial \widetilde I \to t \cdot \partial \widetilde I \=  \{\ ([v], z)\in \widetilde{T_{x_o} M} \ : \ |z| = t \ \}\ , 
\ \  \delta^{(t)}([v], z) = ([v], t z)\ .\eqno(3.6)$$
Since   $\delta^{(t)}$  is a  CR map,  we just need to  show that $\Psi$ maps the real distribution of the 
CR structure of $t \cdot \partial \widetilde I $ onto the real distribution of the CR structure of $S^{(t)}$. 
For this,  let us first observe that  d) implies  $\Psi^* (\tau) = \tau_o$ with $ \tau_o([v], z) \= |z|^2$. Then, let us 
denote by  $J$  the original complex structure of $\widetilde I \subset \widetilde {T_{x_o} M}$ and by  $J'$ the complex structure  
$J' = \Psi^{-1}_*(J)$: The proof reduces now to check that the CR structures induced by $J$ and $J'$ on
$t \cdot \partial \widetilde I = \{ \ x\in \widetilde I\ : \ \tau_o(x) = t^2\ \}$ have the same 
underlying real distribution. But this follows  from Lemma \ref{normaldistribution}
applied to  the pairs $(J, \tau_o)$ and $(J', \tau_o)$.\par 
Claim (iv) is an immediate consequence of the definition of the map $\widetilde f^{(t)}_{([v], z)}$, of properties (a) and (d) of $F$ and the fact that  $\widetilde f^{(t)}_{([v], z)}$ is  holomorphic over the whole $\Delta$, because of boundedness and holomorphicity of $\widetilde f^{(t)}_{([v], z)}$ on $\Delta \setminus \{0\}$,. 
\par
\medskip
To check the last claim,   observe that by definition 
  the map $F: S^{2n-1}Ê\times \Delta \to M$ of 
Thm. 3.4 in \cite{Pt1} extends smoothly to a map $F: S^{2n-1}Ê\times \overline{\Delta} \to \overline{M}$
(see also \S 4 in \cite{Pt3}) and, being the maps  in a) proper,  $F(S^{2n-1} \times \partial \Delta) \subset \partial M$.  From this we infer that 
 $\Psi$ extends smoothly up to the boundary with  $\Psi(\partial \widetilde I) = \partial {\widetilde M}$. Properties (iii) and (iv) when $t = 1$ are obtained  with the same arguments used for  $t < 1$.  \end{pf}
\bigskip
The diffeomorphism  $\Psi: \widetilde I \to \widetilde M$ of the previous proposition 
will be called {\it circular representation of $M$ associated with  $\tau$\/} and 
it can be considered as the lift at blow ups level  of the circular representation 
defined in \cite{Pt1}. When $M$ is a smoothly bounded, strictly convex domain in $\bC^n$,  $\Psi$
coincides with  the circular representation considered in \cite{BDK}.\par
 \bigskip 
 \section{Normal forms for manifolds of circular type}
\bigskip
As in the previous sections,  we will systematically insert the tilde ``\ $\widetilde{\phantom{a}}$\ "  on top of symbols of manifolds,  domains or maps,  whenever 
 we want to indicate a blow up    or a lift of a map at  such blow up. In particular, $\widetilde B^n$ is the blow of the unit ball at the origin. We will also denote by $J_o$ the standard complex structure of $\bC^n$ and of $\widetilde{\bC^n}$ and by $\tau_o: \bC^n \to [0, \infty)$ the standard parabolic exhaustion $\tau_o(z) = |z|^2$.
\par
\subsection{Complex structures of Lempert type and  manifolds in normal form}\hfill\par
\medskip
\begin{definition} \label{Lemperttype} {\rm Let $I \subset \bC^n$ be a complete circular domain with Minkowski function $\mu$ and let  $\pi: \widetilde I\to I$ the  blow up of $I$ at $0$ and $\tilde \mu = \mu \circ \pi$. A complex structure $J$ on $\widetilde I$
is called  {\it of Lempert type\/} if
\begin{itemize}
\item[a)] on any hypersurface $S(c) = \{\ z \in \widetilde I\ :\ \widetilde \mu(z) = c\ \}$,  $0 < c < 1$,  the distribution $\cD$ underlying the CR structure  induced by $J$  coincides with  the   distribution $\cD_o$ underlying the CR structure induced by $J_o$ ; 
\item[b)] the projection $\pi: \widetilde I \to I$ 
induces a  complex manifold structure on $I\setminus\{0\}$, whose charts are smoothly overlapping with 
the charts of the  standard manifold structure of $\bC^n \setminus \{0\}$, i.e.  the projected complex structure  is given by  a smooth tensor field $J$ of type $(1,1)$   on  $I\setminus \{0\}$; 
\item[c)]  the restriction of $J$ on any tangent space of a standard radial disc of $I$ coincides with the standard complex structure $J_o$.
\end{itemize}
Notice that, from a), b) and c), the function $\widetilde \tau \= \mu^2 \circ \pi$  is   plurisubharmonic on   $\widetilde I$ and strictly plurisubharmonic on the complement of the exceptional set.  By Narasimhan's result in  \cite{Na},  if $I$ is endowed with a suitable complex 
manifold structure, this implies that $\widetilde I$ is a proper modification of $I$.
Such complex structure coincides with  the one described in b) on   $I \setminus \{0\}$, but 
 {\it need not\/} to smoothly overlap with the standard complex structure at $0$.\par
We call such complex structure the {\it projected complex structure of Lempert type  on $I$\/} and 
it will be indicated by the associated tensor field $J$, even if  such tensor  is  a smooth tensor w.r.t.  the standard  coordinates  of  $\bC^n$ only at the points of $I \setminus \{0\}$. \par
Two complex structures $J$ and $J'$ of Lempert type   are called {\it Lempert isotopic\/} (or, shortly, {\it L-isotopic\/}) if there exists a smooth family $J_t$, $t \in [0,1]$, of complex structures of Lempert type on $\widetilde I$, such that $J_0 = J$ and $J_1 = J'$.} 
\end{definition}
\begin{theo} \label{maintheorem2}  Any  complex structure $J$  on $\widetilde B^n$,  which  is of Lempert type and L-isotopic to $J_o$,  projects onto a non-standard complex manifold structure $J$ on $B^n$ which makes  $(B^n, J, \tau_o)$ a bounded manifold of  circular type.
\end{theo}
\begin{pf} Since $(B^n,  J_o, \tau_o)$ 
is a domain of  circular type,  we only need to prove  that  conditions  b) and c)  of Definition \ref{circulartype} are still true  after  replacing  $J_o$ with $J$.   
Let $\cZ_o$ and $\cH_o$ be the tangent and normal distributions  of the Monge-Amp\`ere foliation of $\widetilde B^n$ determined by $(J_o, \tau_o)$ and $Z$ the vector field defined  in (2.10).  
Recall that, by remarks after (2.14) and property (a) of Definition \ref{Lemperttype}, 
the distribution $\cH_o$ is $J$-invariant. Recall also that the operator $d^{c'}$,   obtained from 
$d^c$ by replacement of  $J_o$ with  $J$, is
$$d^{c'} = -  \frac{1}{4 \pi} J^{-1} \circ d \circ J\eqno(4.1)$$
where for any $p$-form $\alpha$, $J(\alpha)$ denotes the $p$-form defined by $(J \alpha) (X_1, \dots, X_p) = (-1)^p \alpha(J X_1, \dots, J X_p)$ (see e.g.  \cite{Be}, p. 68). 
 From this we get that, for any vector field $X$   in $\cH_o$,
  $$d d^{c'}\tau_o(Z,  X ) =  \frac{1}{4 \pi} \left\{Z(J X(\tau_o)) -  X (J Z(\tau_o))  - J([X,Z])(\tau_o)\right\} = 0\ .\eqno(4.2)$$
  Here we used the fact that $\cH_o$ is preserved by the flow of $Z$ and that,  by (b) of  Definition \ref{Lemperttype}, 
 $J Z (= J_o Z)$  and $J X$ are  both tangent to  the  hypersurfaces   $\{\ \tau_o = const.\ \}$. \par
Formulae (2.8) and (4.2) show that it suffices  to check  only condition b) of Definition \ref{circulartype}, because this would automatically  imply also  c) of that definition. Now,  $d d^{c'} \tau_o|_{\cZ_o \times \cZ_o} =  d d^{c} \tau_o|_{\cZ_o \times \cZ_o} > 0$, because $(B^n,  J_o, \tau_o)$ satisfies b)  of Definition \ref{circulartype}  and  $J|_{\cZ_o} = J_o|_{\cZ_o}$. On the other hand, 
recall that  $d d^{c'} \tau_o|_{\cH_o \times \cH_o}$ coincides with the Levi form (w.r.t. to $J$) of the  hypersurfaces 
$\{\ \tau_o = const. \ \}$.  Since $\cH_o$ is also the real distribution underlying the CR structure 
 induced by $J_o$ and  the hypersurfaces  $\{\ \tau_o = const. \ \}$ are strongly pseudo-convex (they are  spheres), 
$\cH_o$ is a contact distribution over each such hypersurface (see \S 2.1).  This implies 
that, at any point,  $d d^{c'}\tau_o|_{\cH_o \times \cH_o}$ is a {\it non-degenerate} $J$-Hermitian form.  The same claim is true for all complex structures $J_t$ of an L-isotopy between $J$ and 
$J_o$. A trivial continuity argument implies that  $d d^{c'}\tau_o|_{\cH_o \times \cH_o} > 0$. From this,  (4.8) and (2.7), $d d^{c'} \tau_o > 0$ and $d d^{c'} \log \tau_o \geq  0$, i.e.  b)  of  Definition \ref{circulartype} is  true also when $J_o$ is replaced by $J$.
\end{pf}
 \begin{definition} {\rm
We call  {\it  manifold of circular type in normal form\/} any  bounded manifold of circular type  of the form 
$(B^n, J, \tau_o)$, where  $\tau_o$ is the standard exhaustive function $\tau_o =  |\cdot|^2$ and $J$  is a complex structure of Lempert type that is L-isotopic to $J_o$.\par
Let $(M, J, \tau)$ be a   bounded manifold of circular type. We call {\it normalizing map for $M$ relative to $\tau$ and its center $x_o$\/} any $(J,J')$-biholomorphism $\Phi: \widetilde M \to \widetilde B^n$ between the blow up $(\widetilde M, J)$ at $x_o$ and the blow up $(\widetilde B^n, J')$ at  $0$ of a manifold in normal form $(B^n, J', \tau_o$), so that: 
\begin{itemize}
\item[a)] $\Phi$ induces a diffeomorphism between the exceptional divisors;
\item[b)] $\widetilde \tau =  \widetilde \tau_o\circ \Phi$, where  $\widetilde \tau$ and $\widetilde \tau_o$ are the lifts of  $\tau$ and $\tau_o$ at the blow ups.\end{itemize}}\end{definition}
 \bigskip
 \subsection{Existence and uniqueness of normalizing maps}\hfill\par
\medskip
\begin{theo} \label{firstmaintheorem} Let  $(M, J, \tau)$ be a bounded manifold of circular type,  with  $x_o$ center of $M$ associated with $\tau$ and blow up $\pi: \widetilde M \to M$ at $x_o$. Then,  there
exists at least one normalizing map $\Phi: (\widetilde M, J) \to (\widetilde B, J')$ relative to $\tau$ and 
$x_o$. Moreover,  any  two normalizing maps $\Phi$ and $\hat \Phi$, both relative to $\tau$ and $x_o$,   are equal if and only if    $(\Phi' \circ \Phi^{-1})|_{\pi^{-1}(x_o)} = Id$ and $(\Phi' \circ \Phi^{-1})_*$ induces the identity map on the tangent spaces of the leaves of the Monge-Amp\`ere foliation at the points of $\pi^{-1}(x_o)$.\par
If $(M, J, \tau)$ is a domain of circular type and $\tau$ extends up to the boundary,  then   there exists  a normalizing map which extends smoothly up to the boundary.
\end{theo}
\begin{pf} Fix  a complex basis $(e_0, \dots, e_{n-1})$ for  $T_{x_o} M$ and consider the  unique isomorphism of complex vector spaces $\imath: T_{x_o} M \to \bC^n$ which maps each vector $e_i$ into 
the corresponding vector of the standard basis  $ e^o_i = \imath(e_i)$ of $\bC^n$. In what follows, we constantly identify $T_{x_o} M$  with $\bC^n$
by means of such isomorphism. In particular, use such isomorphism in order to identify the indicatrix $I \subset T_{x_o} M$ associated with $\tau$ with the corresponding  circular domain $I \subset \bC^n$ with Minkowski function $\kappa$.\par
Let $\Psi: \widetilde I \subset \widetilde \bC^n \to \widetilde M$ be the circular representation associated with $\tau$ and consider the complex structure  $J_M = \Psi^{-1}_*(J)$ on $\widetilde I$. Keep in mind that, by (i) of Proposition \ref{circularrepresentation} and definition of $J_M$, the blow up at the origin of  $(I, J_M)$ is precisely  the blow up at the origin of $(I, J_o)$ and that $J_M|_{\bC P^n} = J_o|_{\bC P^n} $. Moreover, 
by construction of $\Psi$ (see formula (3.3)) and property d) in the proof of Proposition \ref{circularrepresentation}, it follows immediately that $\tau \circ \Psi = \kappa^2$ and hence that 
$(I, J_M, \kappa^2)$ is a bounded manifold of circular type. Moreover, 
using  the proof of (iii) in Proposition \ref{circularrepresentation} and by (iv) of the same proposition, it is quite direct to check that $J_M$ satisfies all three conditions for being a complex structure of Lempert type on $I$. We claim that $J_M$ is also L-isotopic to $J_o$. For this, it suffices to consider the diffeomorphisms 
$$\Psi^{(t)}: \widetilde I \to \widetilde M_{< t} = \{\  \tau < t^2\ \}\ ,\qquad \Psi^{(t)}([v], z) = \Psi([v], t z)\ ,\qquad 0 < t < 1$$
and the complex structures $J^{(t)} = \Psi^{(t)-1}_*(J)$ on $\widetilde I$. The same arguments of before show 
that each complex structure $J^{(t)}$ is of Lempert type. 
If we set $J^{(0)}Ê= J_o$ and $J^{(1)}Ê= J_M$, using explicit  coordinate expressions for the maps $\Psi$ and $\Psi^{(t)}$, it can be checked that  $J^{(t)}$, $t \in [0,1]$, is a smooth family of complex structures even at $t = 0$ and $t = 1$, proving that $J_M$ is L-isotopic to $J_o$.\par
 \medskip
 From these remarks, the proof  can be done assuming that $M$  is a smoothly bounded, strongly pseudoconvex circular domain $I \subset \bC^n$, the parabolic exhaustive function $\tau$ is equal to the  $\tau= \mu^2$ where  $\mu$ is the Minkowski function of $I$,  and $J$ is a complex structure of Lempert type on $M = I$, which is L-isotopic to the standard one and so that the blow up of $(I, J)$ at the origin coincides with the blow up of $(I, J_o)$ and with $J|_{\bC P^n} = J_o|_{\bC P^n} $ . \par
 Now, notice that if we replace $J$ by $J_o$, then $(I, J_o, \mu^2)$ remains a  manifold of circular type.  We claim that if $\Phi: (\widetilde I, J_o)  \to (\widetilde B^n, J')$ is a normalizing map for $(I, J_o, \mu^2)$ relative to  $\mu^2$ and center $x_o = 0$, then it is also  a normalizing map also for $(I, J, \mu^2)$.  In fact, it is enough to  consider the complex structure $J'' =  \Phi_*(J_o)$  on $\widetilde B^n$
 and observe that: 
 \begin{itemize}
 \item[1)] $\Phi$ is $(J, J'')$-biholomorphic map by construction; 
 \item[2)]  $\widetilde \tau_o \circ  \Phi = \mu^2\circ \pi$ because $\Phi$ is a normalizing map for $(I, J_o, \mu^2)$; 
 \item[3)] $J''$ is of Lempert type because 
 $J$ is of Lempert type on $I$; 
 \item[4)] $J''$ is L-isotopic to $J'$ (because $J$ is L-isotopic to $J_o$) and 
 $J'$ is L-isotopic to $J_o$; from this it follows  that $J''$ is L-isotopic to $J_o$.
 \end{itemize} 
 So,  if we show the existence of a normalizing map for any smoothly bounded, strongly pseudoconvex circular domain $I \subset \bC^n$, 
 relative to $\tau = \mu^2$ and $x_o = 0$, we automatically prove the existence of 
 normalizing maps for any other manifold of circular type. \par
 \medskip
 This is done using the  lemma that follows. In order to state it, we have 
 to fix some notation. As usual,  the blow-up of $\bC^n$ at the origin is identified  with the tautological 
line bundle $\pi: \widetilde \bC^n = E \longrightarrow \bC P^{n-1}$ and  we set $E_* = E \setminus \{\text{zero section}\} = \bC^n\setminus \{0\}$. We remark that $E_*$   is  
a holomorphic  principal $\bC_*$-bundle over $ \bC P^{n-1}$. \par
Let now  $\mu: \bC^n \to \bR_{\geq 0}$  be the Minkowski function of a smoothly bounded,  complete  circular domain $D \subset \bC^n$ and let $\widetilde \mu = \mu \circ \pi: E_* \to \bR_{\geq 0}$. It is quite direct to realize that  $\tilde \mu^2$  is
the quadratic form of  an Hermitian metric $h^{(\mu)}$ on $\pi: E\to \bC P^{n-1}$ and that  the distribution $\cH$, defined by 
$$\cH_u = \{\ v\in T_u E^*\ :\ d\widetilde \mu_u(v) = d \widetilde \mu_u(J_o v) = 0\ \}\ ,\eqno(4.3)$$
is a connection on the principal $\bC_*$-bundle $E_*$. The associated curvature 2-form $\tilde \omega^{(\mu)}$
is a $\bC_*$-invariant, horizontal 2-form on $E_*$  and  projects down onto a closed,  2-form $\omega^{(\mu)}$ on $\bC P^{n-1}$.  A direct computation shows that, for any  two Minkowski functions $\mu$, $\mu'$, the associated 2-forms
 $\omega^{(\mu)}$ and $\omega^{(\mu')}$  are cohomologous
(see e.g. \cite{PW},  Thm. 3.3). \par
From the definition,  it is clear that, for any hypersurface $S(c) =\{\ \mu = cost.\ \}$, the restriction 
  $\cH|_{S(c)}$  coincides with  the real distribution underlying the  CR structure of $S(c)$. When $D$ is strongly pseudo-convex,  the function $\tau = \mu^2$ is a parabolic exhaustive function for $D$ and hence $\cH|_D$ coincides with   the  normal distribution of the Monge-Amp\`ere foliation of $D$  determined by   $\tau$. Moreover, being each hypersurface   $S(c)$  strongly pseudo-convex, it is simple to check that 
  the associated 2-form $\omega^{(\mu)}$  is a K\"ahler form.\par
  \smallskip
  In the following,  $\mu_o$ is the Minkowski function of $B^n$, i.e. $\mu_o = |\cdot |$,  $\widetilde \mu_o = \mu_o \circ \pi$ and $\cH_o$ is  the corresponding 
  connection on $E_*$ as defined in (4.3).\par
\bigskip
\begin{lem} \label{cruciallemma} Let $D \subset \bC^n$ be a smoothly bounded, strongly 
pseudoconvex circular domain  with Minkowski function $\mu$. Set $\widetilde \mu = \mu \circ \pi$ and 
let $\cH$ be the corresponding connection  on  $E_*$ defined in (4.3). Then, there exists 
a  diffeomorphism $\phi: \widetilde \bC^n \longrightarrow \widetilde \bC^n$  with the following properties:
\begin{itemize}
\item[i)] it is a fiber preserving map for the bundle $\pi: \widetilde \bC^n = E \longrightarrow \bC P^{n-1}$ which is holomorphic on any fiber; 
\item[ii)] $\widetilde \mu_o = \widetilde \mu \circ \phi$; 
\item[iii)] $\phi_*(\cH_o) = \cH$.
\end{itemize}
Moreover, if $\{\mu^{(t)}\ ,\ 0 \leq t \leq 1\}$ is a smooth 1-parameter family of Minkowski functions
of smoothly bounded, strongly pseudoconvex circular domains, 
then it is possible to choose a family of  diffeomorphisms $\phi^{(t)}: \widetilde \bC^n \to \widetilde \bC^n$, satisfying i) -  iii) for $\mu = \mu^{(t)}$,  which   depends smoothly on $t$.
\end{lem}
\begin{pf} Let $\omega_o$ and $\omega$ be the K\"ahler forms on $\bC P^n$ 
determined by the curvatures on $E_*$ of the connections $\cH_o$ and $\cH$, respectively. 
Since $\omega_o$ and $\omega$ are cohomologous, by Moser's theorem \cite{Mo}
there exists a diffeomorphism $\psi: \bC P^{n-1} \to \bC P^{n-1}$ such that 
$$\psi^* \omega = \omega_o\ .\eqno(4.4)$$
Indeed, by the proof of Moser's theorem, there exists a smooth 1-parameter family  of 
diffeomorphisms $\psi_t$, $t \in [0,1]$,  such that $\psi_0 = Id_{\bC P^{n-1}}$, 
$\psi_1 = \psi$. Such family of diffeomorphism is obtained by integrating  a family of  vector fields  $X_t = \dot \psi_t$ satisfying a  particular system of differential  equations. Consider the vector fields  $\hat X_t$ on 
$E_* \simeq \bC ^n \setminus \{0\}$, which are horizontal w.r.t. $\cH_o$ and project onto  the vector fields $X_t$.  Since any space of the distribution $\cH_o$ is tangent to some spheres in  $E_* = \bC^n \setminus \{0\}$, it is possible to integrate such vector fields and obtain another   1-parameter family of diffeomorphisms $\hat \psi_t: E_* \to E_*$ 
so that  $\hat X_t = \dot {\hat {\psi}}_t$.  We  then define  
$$\tilde \psi_t: E \to E\ ,\qquad \tilde \psi_t\left([v], z\right) = \left\{\begin{matrix} (\psi_t([v]), \hat \psi_t(z))&\text{if}\ z \neq 0\\
\phantom{a}&\phantom{a}\\
(\psi_t([v]), 0)&\text{if}\ z = 0
\end{matrix}
\right.\ ,$$
which  can be easily checked  to be a family  of diffeomorphisms and we set 
$\widetilde \psi \= \hat \psi_1$. By construction, $\widetilde \psi: E \to E$ commutes with the action of $\bC$ on $E$ and 
restricts to $\psi$ on $\bC P^{n-1}$.
Consider now a  diffeomorphism $\phi: E \to E$ of the form 
$$\phi([v], z) \= \left(\psi([v]),  e^{i \lambda([v])} \frac{|z| \cdot \hat \psi(z) }{\mu(\hat \psi(z))}\right)\eqno(4.5)$$
where we denote by $\lambda: \bC P^{n-1} \to \bR$ a smooth function to be fixed later.
It is clear that  a map  of the form (4.5) satisfies  (i). But it satisfies also (ii), since 
$$\widetilde \mu\left(\phi([v], z)\right) = \mu\left(e^{i \lambda([v])} \frac{|z| \cdot \hat \psi(z) }{\mu(\hat \psi(z))}\right) = |z| = \widetilde \mu_o([v], z)\ ,$$
and we claim that  there exists a $\lambda$ so that 
it satisfies also  (iii). To check this, observe that, since  $\phi|_{E_*}$ commutes with the action of 
$\bC_*$,  it maps the connection $\cH_o$ into a connection $\cH'$ on $E_*$ whose curvature 
2-forms   projects down on $\bC P^{n-1}$ onto the $\bC$-valued 2-form 
$$ \omega' = \psi_*\omega_o = \omega\ .\eqno(4.6)$$
Moreover, by (ii),  all spaces of $\cH'$ are tangent to the hypersurfaces $\{\ \mu = cost.\ \}$.\par
By standard arguments of theory of connections,  we obtain that, for any $u \in E_*$
there exists a linear map
$$\nu_u: T_{\pi(u)} \bC P^{n-1} \longrightarrow \bR \eqno(4.7)$$
so that any space $\cH'_u$   is of the form
$$\cH'_u = \left\{\ w = v +  \nu_u(\pi_*(v)) \cdot \left(i \left.\frac{\partial}{\partial \zeta}\right|_u - i \left.\frac{\partial}{\partial \bar  \zeta} \right|_u \right)\ \text{for some}\  v \in \cH_u\ \right\}\eqno(4.8)$$
(here we denoted by $\zeta$ the standard coordinate of $\bC$ and by $\frac{\partial}{\partial \zeta}$ the 
vertical holomorphic vector field on $E$ induced by the holomorphic action of $\bC$ on $E$). 
Equivalently, we may say that if $\varpi: TE_* \to \bC_*$ is the connection form for $\cH$, then 
the connection form $\varpi'$ of $\cH'$ is 
$$\varpi' = \varpi -  i \nu\ .\eqno(4.9)$$
By the invariance of $\cH'$ and $\cH$ under the $\bC_*$-action,   the linear map $\nu_u$ depends only on the point $x = \pi(u) \in \bC P^{n-1}$   and  we may consider $\nu$ as a 1-form on $\bC P^n$. Computing the curvature,  we get from (4.9) and (4.6) that 
$$d\nu = i(\omega' - \omega) =  0\ .\eqno(4.10)$$
Since $H^1(\bC P^{n-1}) = 0$, there exists a smooth function $\widetilde \lambda: \bC^n \to \bR$ such that 
$\nu = d \widetilde \lambda$. Now, let us replace the function $\lambda$ in (4.5) with the function 
$\lambda - \widetilde \lambda\circ \psi^{-1}$. By construction, in (4.8) the function $\nu_u$ 
has to be replaced by the function $\widetilde \nu_u = \nu_u - \left.d\widetilde\lambda\right|_{\pi(u)} = 0$
and   the  new map $\phi$ satisfies (iii).\par
It remains to prove the final part of the statement. First of all, notice that by the first part of the proof, 
the function  $\lambda: \bC P^{n-1} \to \bR$ is uniquely determined (up to a constant) by the diffeomorphism $\psi: \bC P^{n-1} \to \bC P^{n-1}$ which satisfies (4.4). By choosing some suitable normalization condition for 
$\lambda$,  we may assume that the map $\phi$   is uniquely determined 
by   $\psi$. 
If $\mu^{(t)}$ is a smooth family of Minkowski functions of strongly pseudoconvex, 
complete circular domains, 
then also the  corresponding K\"ahler  forms $\omega^{(t)}$  are smoothly depending 
on $t$ (it suffices to see the explicit 
expression of $\omega^{(t)}$ in term of $\mu^{(t)}$ -  see e.g. \cite{PW}, p.27). Now, 
the  proof of Moser's theorem in \cite{Mo} shows that there exists a smooth family of diffeomorphisms $\psi^{(t)}$, $t \in [0,1]$, each of them  satisfying $\psi^{(t)*} \omega^{(t)} = \omega_o$.  This automatically 
implies the existence of a smooth family $\phi^{(t)}$  satisfying (i) - (iii) for any $t$.
\end{pf}
We may now conclude the proof of the theorem.  Given a smoothly bounded, strongly pseudoconvex circular domain $I \subset \bC^n$ with Minkowski function $\mu$, let 
$\Phi = \phi^{-1}|_{\widetilde I}: \widetilde I \to \widetilde B^n$, where $\phi$ is diffeomorphism   of the previous lemma,  and let $J' = \Phi_*(J_o)$. From (i) - (iii) of the lemma, it follows immediately that $J$ is of Lempert type and that $\widetilde \tau = \widetilde \tau_o \circ \Phi$ with  $\tau = \mu^2$. Moreover,  we may consider the smooth 1-parameter family of Minkowski functions $\mu^{(t)} = (1-t) \mu + t \mu_o$ with  $t \in [0,1]$, a corresponding family of diffeomorphisms $\Phi_t = \phi_t^{-1}$ with $\phi_t$ associated by the lemma with $\mu^{(t)}$ and smoothly depending on $t$, and the 1-parameter family of complex functions $J_t \= \Phi_t{}_*(J_o)$. By construction,  $J_t$ is an L-isotopy between $J'$ and $J_o$ and we may conclude that $\Phi$ is a normalizing map for $(I, J_o, \tau = \mu^2)$,  as needed.\par
\smallskip
Assume now that $\Phi, \hat \Phi: \widetilde M \to \widetilde B^n$ are two normalizing maps  so that $(\Phi'{}^{-1} \circ \Phi)|_{\pi^{-1}} = Id_{\pi^{-1}}$ or, equivalently, that $(\Phi'\circ \Phi^{-1})|_{\bC P^{n-1}} = Id|_{\bC P^{n-1}}$.  Now,  $\Phi$ and $\Phi'$ map the leaves of the Monge-Amp\`ere foliation of $(M, J, \tau)$ into the 
leaves of the Monge-Amp\`ere foliation of $(B^n, \Phi_*(J), \tau_o)$ and of $(B^n, \Phi'_*(J), \tau_o)$, which are in both cases the images of the standard radial discs.  This means that 
$\Phi'\circ \Phi^{-1}$ maps biholomorphically any standard radial disc into itself, since 
  $(\Phi'\circ \Phi^{-1})|_{\bC P^n} = Id|_{\bC P^{n-1}}$. If in addition $(\Phi'{}^{-1} \circ \Phi)_*$ induces the identity map on any tangent space of a leaf of the Monge-Amp\`ere foliation at the points of $\pi^{-1}(x_o)$, we get that   $\Phi'\circ \Phi^{-1}$ maps any radial disc into itself by  a biholomorphism which  fixes the origin and  with derivative  equal to $1$ at $0$.  By Schwarz lemma,  $\Phi'\circ \Phi^{-1} = Id$ on any radial disc and hence on the whole $B^n$.\par
\smallskip
It remains to check the smooth extendibility up the boundary of $\Phi$ if $\tau$ is smoothly extendible. But this  is a consequence of the  fact that  $\Phi$ is obtained by composing the inverse of the circular representation (which is smoothly extendible to the boundary because of Proposition \ref{circularrepresentation}) and the diffeomorphism between  $\widetilde I \subset \widetilde T_{x_o} M = \widetilde \bC^n$ and $\widetilde B^n$,  which is given in Lemma \ref{cruciallemma}  and which is trivially  smooth  up to the boundary.
\end{pf}
\smallskip
We remark that one can prove a stronger statement 
about the uniqueness of normalizing maps. We will come back on this topic in \S 6.
\bigskip
\section{Bland and Duchamp's invariants}
\bigskip
In this section, we consider only manifolds of circular type in normal form, i.e. 
of the form  $(B^n, J, \tau_o)$ with $\tau_o = |\cdot |^2$ and $J$  complex structure 
on $\widetilde B^n$ of Lempert type and L-isotopic to the standard one.
The distribution in $\widetilde B^n$, 
which is  normal to the standard Monge-Amp\`ere foliation,  will be denoted $\cH$. We recall that, for  any sphere $S(c) = \{\ \tau_o = c^2\ \}$, the restriction $\cH|_{S(c)}$ 
is  the distribution underlying the standard CR structure of $S(c)$  and that,  for any $x \in \widetilde B^n$, 
$$T_x \widetilde B^n = \cZ_{x} \oplus \cH_x$$
where $\cZ$ is the distribution tangent to the standard radial disc. We also denote by $Z$ the vector field defined in (2.10):  in the standard coordinates of $\bC^n$, the corresponding holomorphic and anti-holomorphic parts of $Z$ are 
$$Z^{1,0} = \frac{1}{2}\left(Z - i J_o Z\right) = z^i \frac{\partial }{\partial z^i}\ ,\qquad Z^{0,1} =  \bar z^i \frac{\partial }{\partial \bar z^i}\ .\eqno(5.1)$$
Recall also that any complex structure $J$ of Lempert type is uniquely determined by its action on the vector fields  on  $\cH$, since 
the action on the vector fields in $\cZ$ is the same of the standard complex structure $J_o$.\par
\smallskip
Let $H^{1,0}$ and $H^{0,1} = \overline{H^{1,0}}$ be the $J_o$-holomorphic and $J_o$-anti-holomorphic subbundles 
of $\cH^\bC$. For any other complex structure $J$ of Lempert type,  we denote 
the  corresponding   $J$-holomorphic  and $J$-antiholomorphic subbundles of $\cH^\bC$ by $H^{1,0}_J$ and $H^{0,1}_J = \overline{H^{1,0}_J}$.\par
\begin{definition} {\rm Let  $J$ be a complex structure on $\widetilde B^n$ of Lempert type.  We call  {\it deformation tensor associated with $J$\/} 
any smooth section
$$\phi: \widetilde B^n \to \bigcup_{x \in \widetilde B^n} \operatorname{Hom}(H^{0,1}, H^{1,0}) = H^{0,1}{}^*\otimes H^{1,0}$$  
so that $H^{0,1}$ can be expressed as  
$$H^{0,1}_J|_x = \{\  w + \phi_x(w)\ ,\ w \in H^{0,1}|_x\ \}\ \qquad \text{for any}\ x \in \widetilde B^n\ .\eqno(5.2)$$
\/} 
\end{definition}
\medskip
Notice that, a priori, not any complex structure of Lempert type has an associated deformation tensor. However  if  $J$  has an associated deformation tensor, 
then any sufficiently small deformation $J'$ of $J$, which is also of Lempert type, 
has an associated deformation tensor. Indeed, we will shortly see that  {\it any complex structure 
$J$ of Lempert type and  L-isotopic to $J_o$  admits an associated  deformation tensor\/}. \par
\medskip
We now want to exhibit some differential equations which characterize the deformation tensors. In order to do  this, we first 
need to recall the definition of  two important operators on the tensor fields in  $H^{0,1}{}^*\otimes H^{1,0}$. \par
We recall that $\widetilde B^n \subset \widetilde \bC^n$
is a holomorphic bundle $\hat \pi: \widetilde B^n \to \bC P^{n-1}$, with fibers given by the radial discs (if endowed with the standard complex structure). Since $\cH$ is a connection 
in such a bundle, the holomorphic and anti-holomorphic distributions  are generated  by vector fields $X^{1,0} \in H^{1,0}$ and $Y^{0,1} \in \cH^{0,1}$ that can be locally chosen 
so that $\hat \pi_*\left([X^{1,0}, Y^{0,1}]\right) = [\hat \pi_*(X^{1,0}),\hat \pi_*(Y^{0,1})] = 0$. Let us call such vector fields {\it holomorphic\/} and {\it anti-holomorphic\/} vector fields 
of the distribution $\cH^\bC$, respectively. 
It can be easily checked that if $\phi$ is a deformation tensor associated with a complex structure of Lempert type, then for any two anti-holomorphic vector fields
$X, Y \in H^{0,1}$
$$[X, \phi(Y)] \in H^{1,0} + \cZ^\bC\ .\eqno(5.3)$$
Hence, if we denote by $( \cdot )_{\cH^\bC}$ the projection onto the distribution $\cH^{\bC}$, we surely have that 
 $[X, \phi(Y)]_{\cH^\bC} \in H^{1,0}$ for any pair of anti-holomorphic vector fields.
Now,  consider the following two operators (see \cite{KM}):
$$\bar \partial_b: H^{0,1}{}^*\otimes  H^{1,0} \to \Lambda^2 H^{0,1}{}^* \otimes  H^{1,0}\ ,\ $$
$$
\bar \partial_b \alpha(X,Y) \= [X,\alpha(Y)]_{\cH^\bC} -    [Y,\alpha(X)]_{\cH^\bC} -
\alpha([X,Y]) 
\ ,\eqno(5.4)$$
and 
$$[ \cdot, \cdot ] : \left(H^{0,1}{}^*\otimes  H^{1,0}\right) \times \left(H^{0,1}{}^*\otimes  H^{1,0} \right) \longrightarrow \Lambda^2 H^{0,1}{}^* \otimes  H^{1,0}\ ,$$
$$[\alpha, \beta](X,Y) \= \frac{1}{2} \left([\alpha(X), \beta(Y)] - [\alpha(Y), \beta(X)]\right)\ \eqno(5.5)$$
for any pair of anti-holomorphic vector fields $X$, $Y$ in  $H^{0,1}$. \par
\begin{prop} \label{invariants} Let $J$ be a  complex structure on $\widetilde B^n$ of Lempert type that admits an associated deformation tensor $\phi$.
Then:
\begin{itemize}
\item[i)] $dd^c\tau_o(\phi(X), Y) + dd^c \tau(X, \phi(Y)) = 0$ for anti-holomorphic  $X,Y \in H^{0,1}$; 
\item[ii)] $\bar \partial_b \phi + \frac{1}{2}[\phi, \phi] = 0$; 
\item[iii)]  $\cL_{Z^{0,1}}(\phi) = 0$. 
\end{itemize}
Conversely, any tensor field $\phi \in H^{0,1}{}^*\otimes H^{1,0}$ that satisfies (i) - (iii) is the deformation tensor
of a complex structure of Lempert type. \par
In addition, a complex structure $J$ of Lempert type, associated with a deformation tensor $\phi$,
is so that $(B^n, J , \tau_o)$ is a manifold of circular type if and only if 
\begin{itemize}
 \item[iv)] $ dd^c \tau_o(\phi(X), \overline{\phi(X)}) < dd^c \tau_o(\bar X, X)$   for any  $0 \neq X \in H^{0,1}$.
 \end{itemize}
\end{prop}
\begin{pf} First of all, recall  that by the $J_o$-invariance of the 2-form $d d^c \tau_o$,  for any two  vector fields in $H^{0,1} \oplus \cZ^{0,1}$ or 
in $H^{1,0} \oplus \cZ^{1,0}$,
$$d d^c\tau_o(W, W') = d d^c\tau_o(J_o W, J_o W') = - d d^c\tau_o(W, W') =  0\ .$$
So, from the proof of Theorem \ref{maintheorem2}, the reader can check  that a complex structure $J$ of Lempert type is so that $(B^n, J, \tau_o)$ is a manifold of circular type, if and only if for any  $0 \neq X \in H^{0,1}$
$$dd^c \tau_o(\bar X + \overline{\phi(X)}, X + \phi(X)) = dd^c \tau_o(\bar X,  X) + dd^c \tau_o( \overline{\phi(X)}, \phi(X))> 0\ .$$
This proves (iv). 
For  checking the necessity and sufficiency of (i) - (iii), we only need to show that those properties are necessary and sufficient condition
for the integrability of the unique  almost complex structure $J$, which coincides with $J_o$ on the radial discs, leaves the distribution  $\cH$ 
invariant and have an associated   anti-holomorphic distribution  $H^{0,1}_J$ which is  as in (5.2).
Such almost complex structure   $J$ is integrable if and only if for any anti-holomorphic vector fields  $X, Y \in H^{0,1}$ 
one has 
$$[X + \phi(X), Y + \phi(Y)]Ê\in \cZ^{0,1} + H^{0,1}_J\ ,\qquad [Z^{0,1}, X + \phi(X)] \in \cZ^{0,1} + H^{0,1}_J\ .\eqno(5.6)$$
But conditions  (5.6) are satisfied if and only if 
$$ [X + \phi(X), Y + \phi(Y)]_{\cH^\bC}  =  [X,Y] + \phi([X,Y]) \ \Leftrightarrow\ \bar \partial_b \phi(X,Y) + \frac{1}{2}[\phi, \phi](X,Y) = 0\ ,\eqno(5.7)$$
$$  [X + \phi(X), Y + \phi(Y)]_{\cZ^\bC} = 0 \ \Leftrightarrow \ d d^c \tau_o(X + \phi(X), Y + \phi(Y)) = 0\ ,\eqno(5.8)$$
$$[Z^{0,1}, X + \phi(X)] = [Z^{0,1}, X] + \phi([Z^{0,1}, X] ) \ \Leftrightarrow\ \cL_{Z^{0,1}} \phi (X) = 0\ \eqno(5.9)$$
for any anti-holomorphic  $X$, $Y \in H^{0,1}$, i.e. if and only if (i) - (iii) are true.
\end{pf}
\bigskip
Let $J$ be a complex structure of Lempert type and $J_t$, $t \in [0,1]$, an  L-isotopy between $J$ and $J_o$. By the previous remark, the set of $t$'s, for which $J_t$ has an associated  deformation tensor 
is open, while (iv) of the previous lemma implies that it is also closed. From this, we conclude that 
also $J = J_1$ has a deformation tensor and hence that  {\it there is a natural injective map between  the class of manifolds in normal form $(B^n, J, \tau_o)$  and the class of tensor fields  $\phi \in H^{0,1}{}^* \otimes H^{1,0}$ on $\widetilde B^n$ which satisfy (i) - (iv) of Proposition \ref{invariants} \/}. \par
The correspondence between normal forms and deformation tensors satisfying (i) - (iv) is a priori only injective, not surjective.  However,  for any deformation tensor  satisfying   (i) - (iv), the associated complex structure $J$ defines a manifold of circular type and hence  there exists some  normalizing map $\Phi: \widetilde B^n \to \widetilde B^n$ for which $\hat J = \Phi_*(J)$ is in normal form and whose associated deformation tensor is $\hat \phi = \Phi_*(\phi)$. In other words, we may say that {\it any deformation tensor satisfying  (i) - (iv) is, up to a diffeomorphism,  the  deformation tensor of some normal form\/}.\par
\medskip
 Proposition \ref{invariants} (iii) has also the following consequence. Consider $n-1$ holomorphic vector fields $(e_1, \dots, e_{n-1})$, defined on some open subset of $\bC P^{n-1} \subset \widetilde B^n$ and linearly independent at all points where they are defined. 
 Extend  them to $\widetilde B^n$  as $Z$ and $JZ$ invariant vector fields on taking values in  $H^{1,0}$. 
 Let also $(e^1, \dots, e^{n-1}$, $Z^{1,0}{}^*)$ the holomorphic  field of $(1,0)$-forms, which is dual to the frame field $( 
 e_1, \dots, e_{n-1}$, $Z^{1,0})$. Then any tensor field $\phi \in H^{0, 1}{}^* \otimes H^{1,0}$ is of the form 
 $\phi = \sum_{a,b = 1}^{n-1}\phi^a_b \overline{e^b}\otimes e_a$ and satisfies  Proposition \ref{invariants} (iii) if and only if the restrictions of the functions $\phi^a_b$ on the radial discs are  holomorphic. 
 In particular, using a system of coordinates $(v^1, \dots, v^{n-1}, \zeta)$ for $\widetilde B^n$ as in (2.5),  we have that $\phi$ satisfies (iii) if and only if it  is of the form 
 $$\phi = \sum_{j = 0}^\infty \phi_j \zeta^j \ ,\qquad \phi_j =  \sum_{a,b = 1}^{n-1}\phi^a_{bj} \overline{e^b}\otimes e_a$$ 
 where $\phi^a_{bj} = \phi^a_{bj}(v^1, \dots, v^{n-1})$ are the coefficients of the series expansion in powers  of $\zeta$ of the functions 
 $\phi^a_{b}(v^1, \dots, v^{n-1}, \zeta)$.  It can be checked that the deformation tensors $\phi^{(k)} \= \phi_k \zeta^k \in H^{0,1}{}^* \otimes H^{1,0}$ are independent on the choice of the coordinates and of the frame field $(e_1, \dots, e_{n-1})$. Moreover,  $\phi$ satisfies (i) and (iii) of Proposition \ref{invariants}  if and only 
 each tensor field $\phi^{(i)}$ satisfies (i) and that the following  equations:
 $$\bar \partial_b \phi^{(k)} + \frac{1}{2} \sum_{i + j = k} [\phi^{(i)}, \phi^{(j)}] = 0\qquad \text{for any} \ 0 \leq k < \infty\ . \eqno(5.10)$$
  Summarizing, we have the following theorem, which can be considered as an  extension  to  arbitrary manifolds of circular type some of the main results in \cite{Bl, BD1} (see next remark).\par
 \medskip
 \begin{theo} \label{BlandandDuchamp}ÊLet   $D$ be a manifold  of circular type  in normal form, i.e. $D = (B^n, J, \tau_o)$, with $J$ complex structure of Lempert type and L-isotopic to $J_o$.  Then   $J$ is uniquely
 determined by an associated  sequence of deformation tensors $ \phi^{(k)} \in H^{0,1}{}^* \otimes H^{1,0}$, $0 \leq k < \infty$,  which   satisfy (5.10) and  (i) of Proposition \ref{invariants} for any $k$, and 
 with the series $\phi = \sum_k \phi^{(k)}$ uniformly  converging on compacta.
 \end{theo}
\medskip
In   \cite{Bl, BD1}, Bland and Duchamp considered small deformations of the standard CR structure of the 
 $S^{2n-1}$ and proved that,  for $n >1$, any such CR structure  is embeddable in $\bC^{n-1}$ as boundary of a domain which is biholomorphic to a domain in normal form 
$(B^n, J, \tau_o)$.  In particular, they associated with any small deformation of the CR structure of $S^{2n-1}$ a sequence of tensors,  which correspond to the restrictions to $S^{2n-1}$ of the  tensors $\phi^{(k)}$ appearing in   Theorem \ref{BlandandDuchamp}
\footnote{Be aware  that  our deformation tensor is  minus the deformation tensor considered in  \cite{Bl, BD1}.}.  It is therefore natural to name such sequence of  deformation tensors $\phi^{(k)}$  the {\it Bland and Duchamp invariants of $(B^n, J, \tau_o)$\/}.\par
\medskipÊ 
We conclude with the following concept, which will turn out to be  quite useful in the applications of the next section.\par
\begin{definition} {\rm Let $(D, J, \tau)$ be a domain of circular type.  We say that  $D$ is {\it stable\/}   if  it admits a parabolic exhaustive function $\tau$, whose associated  circular representation 
$\Psi: \widetilde I \to D$  extends smoothly up to the boundary inducing  a 
 diffeomorphism   between  $\partial \widetilde I$ and $\partial D$. \par
Unless differently stated, for any  given stable domain $D$, we will call {\it parabolic exhaustive  functions of $D$\/} only those, whose associated circular representation  satisfies the above condition.} 
\end{definition}
\smallskip
 We also recall that,
by  Lempert's  and  the first author's result (see e.g. \cite{Le, Pt}),  the class of  stable domains of circular type naturally includes the  smoothly bounded, strictly convex domains of $\bC^n$ and  the 
smoothly bounded, strongly pseudoconvex circular domains.  Indeed,  by  Thm. 4.4 in \cite{Pt1}, the complete circular domains may be characterized as the  unique (stable) domains  of circular type domains whose  Monge-Amp\`ere foliation is holomorphic. \par
  Moreover, 
the  ``stability'' property   is invariant under biholomorphisms between domains of circular type. In fact:
\par
\begin{lem} Let $(D, J, \tau)$  and $(\hat D, \hat J, \hat \tau)$ be two biholomorphic  domains  of circular type. Then $D$ is   stable if and only if   $\hat D$ is stable. Furthermore  they  have a 
normal form $(B^n, J, \tau_o)$, with  a complex structure $J$ which is smoothly extendible up to the boundary and makes  $\overline{B^n}$ a stable domain of circular type.
\end{lem}
\begin{pf} The first claim follows from  the fact  that any biholomorphism $f: D \to \hat D$ between stable domains extends smoothly up to the boundary. This  property can be checked   using  the  local regularity results  of  Berteloot (\cite{Ber},  Prop. 3), which imply that  any such  $f$  admits an
 H\"older continuous extension up the boundary.  In fact, from H\"older boundary regularity,   the standard arguments of Lempert's proof of Fefferman theorem (see \cite{Le, Tu})  imply    that $f$ extends smoothly up to the boundary.\par
The last claim is a consequence of the proof of  Theorem \ref{firstmaintheorem}. In fact, if
$D$ is stable, we may construct a normalizing map which is smooth up to the boundary and induces a diffeomorphism between $\partial D$ and $S^{2n-1} = \partial B^n$. In particular, the complex structure of $\widetilde D$, which extends up to $\partial D$, is  mapped onto a complex structure on $\widetilde B^n$, which extends smoothly up to $\partial B^n$ (and hence also to a small neighborhood of $\overline{B^n}$). This implies that the associated normal form $(B^n, J ,\tau_o)$ is a stable domain of circular type since the circular representation of  $(B^n, J ,\tau_o)$ coincides with the one 
of $(B^n, J_o ,\tau_o)$.
\end{pf}
\medskip
From the previous lemma, it is clear that the normal forms of stable domains of circular type 
are associated with sequences of Bland and Duchamp's invariants $\{\phi^{(k)}\}$ which converge uniformly on the closure $\overline{\widetilde B^n}$.\par
\bigskip
\section{Miscellaneous results}
\bigskip
\subsection{The geometrical meaning of the Bland and Duchamp invariant  $\phi^{(0)}$}\hfill\par
\label{BD-invariant0}
\medskip
Let $(B^n, J, \tau_o)$ be a manifold in normal form, $\phi^{(k)}$ the associated Bland and Duchamp invariants, $I \subset T_0 B^n \simeq \bC^n$ the indicatrix at $0$ and $\Psi: \widetilde I \to \widetilde B^n$ the 
circular representation.  The identification $T_0 B^n \simeq \bC^n$ is done so that we may assume $J|_{T_0 B^n} = J_o$.  
The following proposition collects a few properties which give a clear indication of  the information carried by  the Bland and Duchamp invariant $\phi^{(0)}$.\par
\begin{prop}\label{propositioncircular}\hfill \par
\begin{itemize}
\item[a)] The pull-backed complex structure   $J' = \Psi^*(J)$  on $\widetilde I$ is  of Lempert type.
\item[b)]  The tensor field $\phi - \phi^{(0)} = \sum_{k \geq 1}^\infty \phi^{(k)}$ is  identically vanishing   if and only if  the circular representation  is a biholomorphism between $I$   and $(B^n, J, \tau_o)$,   i.e. if and only if $(B^n, J)$ is biholomorphic to a circular domain. 
\item[c)]  The invariant $\phi^{(0)}$ is   always the deformation tensor of a manifold in normal form $(B^n, J^{(0)}, \tau_o)$, more precisely,  of a normal form of  the indicatrix $I $.\end{itemize}
\end{prop}
\begin{pf} a) is a direct consequence of definitions and Proposition \ref{circularrepresentation}. For b), we remark that $\phi^{(k)} = 0$ for all  $k \geq 1$ if and only if the projection $\pi: \widetilde B^n \to \bC P^{n-1}$ is a $J$-holomorphic, i.e. if and only if the Monge-Amp\`ere foliation of $(B^n, J, \tau_o)$ is holomorphic. Then (b) follows from \cite{Pt1}, Prop. 3.4. \par
For c), notice that   $\phi^{(0)}$ satisfies (i) - (iv)  of Proposition \ref{invariants} and defines a complex structure which is L-isotopic to $J_o$,  because $\phi$ does it. So, by the remarks after Proposition \ref{invariants}, the first claim follows immediately. Now, consider the circular representation  $\Psi: \widetilde I \to \widetilde B^n$. It is straightforward to realize that $J^{(0)}|_{\bC P^n} = \Psi_*(J_o|_{\bC P^n})$ and hence, by invariance along the leaves of the Monge-Amp\`ere  foliations,
$J^{(0)} = \Psi_*(J_o)$ on $\widetilde I$. This implies that the corresponding projected structures on $I$ and on $B^n$ are biholomorphic.
 \end{pf}
 \bigskip
 \subsection{The parameterization of normalizing maps and the automorphism group of a manifold of circular type}\hfill\par
\subsubsection{Special frames of a manifold of circular type}\hfill\par
\begin{definition} \label{specialframe} {\rm Let $\tau$ be a parabolic exhaustion function  for $M$ and 
$x_o$ and $I_{x_o} \subset T_{x_o} M \simeq \bC^n$ the corresponding center and indicatrix. We call {\it special frame at $x_o$ associated with $\tau$\/} a complex basis $(e_0, e_1, \dots, e_{n-1})$ for 
$T_{x_o} M$
defined as follows:
\begin{itemize}
\item[i)] $e_0 \in \partial I_{x_o}$ (i.e. $\kappa(e_0) = 1$, where $\kappa$ is the Minkowski function of $I_{x_o}$);
\item[ii)] $(e_1, \dots, e_{n-1})$ is a unitary basis w.r.t. $d d^c \kappa^2$ for the holomorphic tangent space   $\cD^{1,0}_{e_0} \subset T_{e_0} \partial I_{x_o}$ 
of the CR structure of $ \partial I_{x_o} \subset T_{x_o} M \simeq \bC^n$.
\end{itemize}
Recall that if $D$ is a {\it domain\/} of circular type,  for any center $x_o$ there is a unique  parabolic exhaustion function $\tau: D \to [0,1)$, smoothly extendible at the boundary,  for  which the center is exactly  $x_o$ (see Lemma \ref{uniqueness}). For this reason, for any such domain the following set is well defined
$$P = \bigcup_{x_o \text{is a center}} P_{x_o}\ ,\qquad\text{where}\qquad  P_{x_o} = \{\ \text{special frames at} \ x_o\ \}\ .$$
We will call it  the {\it pseudo-bundle of special frames of $D$}. We will also denote by $\pi: P \to D$ the map which associates to any special frame  the base point  and ${\mathfrak C}(D) = \pi(P) \subset D$ will denote the set of centers of $D$.}
\end{definition}
\medskip
It should be observed that, relatively to  the  biholomorphisms that are  smooth up to the boundary (in case of stable domains, any biholomorphism is in such a class),  the pseudo-bundle $P$ is a biholomorphic  invariant of   $D$. 
In case  $\mathfrak C(D)$ is discrete,  the pseudo-bundle $P$ is a bundle over such a set. Any   fiber  $P_{x_o} = \pi^{-1}(x_o)$ has a structure of  $\UU_{n-1}$-principal bundle with basis equal  to the boundary of the indicatrix $\partial I_{x_o}$. In case $D$ is a smoothly bounded, strictly convex 
domain in $\bC^n$ (and hence  $D = \mathfrak C(D)$), $P$  is a bundle over $D$. More precisely,  it coincides with   {\it unitary frame bundle of the complex Finsler metric\/}  given by the Kobayashi metric of $D$ (\cite{Sp}). \par
Notice also that, via Gram-Schmidt ortonormalization, any special frame of a manifold in normal form is uniquely associated with a basis which is unitary w.r.t. the standard Poincar\`e-Bergman metric of the unit ball   $B^n \subset \bC^n$. Using such correspondence, it is possible to embed the pseudo-bundle $P$ into the bundle $\operatorname{U}_n(B^n)$ of the unitary frames of $B^n$. Since $\Aut(B_n, J_o)$ acts transitively and freely on $\operatorname{U}_n(B^n)$,  $ \operatorname{U}_n(B^n)$ can be identified  with $\Aut(B_n, J_o)$ and the previous immersion  $P\hookrightarrow  \operatorname{U}_n(B^n) $ can be considered as an immersion $P\hookrightarrow \Aut(B_n, J_o)$. In case $\mathfrak C(D) = D$, such immersion is actually a diffeomorphism (see e.g. \cite{Sp, ST}).\par
\bigskip
\subsubsection{Parameterization of normalizing maps by means of special frames}\hfill\par
\medskip
Consider a manifold of circular type $(M, J, \tau)$ of dimension $n$ and denote by $\widetilde {\cN(M)}$ the class of all normalizing maps $\Phi: \widetilde M \to \widetilde B^n$ between a blow up of $M$ at a center $x_o$ and the blow up of $B^n$ at the origin. A priori, $\widetilde {\cN(M)}$ 
is a very large class of maps: In fact, for any  $\Phi \in \widetilde {\cN(M)}$, if $\psi: \widetilde B^n \to \widetilde B^n$ is a diffeomorphism, 
which is  a fiber bundle automorphism of $\pi: \widetilde B^n \to \bC P^{n-1}$ and preserves the complex structure on the fibers, the composition  $\psi\circ \Phi$ is still in $ \widetilde {\cN(M)}$. This means that, for any given $\tau$ and $x_o$, the class  $\widetilde {\cN(M)}$ contains a family   of maps of  cardinality larger or equal to   the cardinality of  $\hbox{\it Diff}(\bC P^{n-1})$. \par
On the other hand,  any normalizing map $\Phi: \widetilde M \to \widetilde B^n$ induces a complex structure on $\widetilde B^n$, which in turn projects down onto a complex manifold structure on $B^n$ (see Definition \ref{Lemperttype}).  In general, the charts of two distinct complex manifold structures on $B^n$ of this kind do not smoothly overlap, i.e.   {\it belong to two distinct   real manifolds structures\/}, even if they are surely diffeomorphic and they coincide with  the standard manifold structure of $\bC^n$ when  restricted to $B^n \setminus \{0\}$ .\par
For this reason, in the following, given a manifold $M$ of circular type,  we fix   {\it once and for all\/}  one of  such real manifold structures and we will denote by $\cN(M) \subset \widetilde \cN(M)$  the class of  normalizing maps which induce on $B^n$ that   real manifold structure. In other words,  a normalizing map $\Phi : \widetilde M \to \widetilde B^n$ belongs to  $\cN(M)$ if and only if the  corresponding  projected map $\phi: M \to B^n$ is a diffeomorphism w.r.t. to the real manifold structure of $M$ and the fixed manifold structure on $B^n$.  Finally, for any parabolic exhaustive function $\tau$ on $M$ and corresponding center $x_o$, we denote by $\cN(M, \tau, x_o)$ the subclass of  $\cN(M)$ consisting of normalizing maps which are associated with   $\tau$ and   $x_o$.  \par
\medskip
\begin{lem}\label{lemmanormalmaps} Let $\Phi, \Phi' \in \cN(M, \tau, x_o)$ and denote by $\phi, \phi': M \to B^n$
the corresponding projected maps. Then $\Phi = \Phi'$ if and only if 
$\left.\left(\phi'\circ \phi^{-1}\right)_*\right|_{T_0 B^n} = Id$.\end{lem}
\begin{pf} The necessity is immediate. For the sufficiency, notice that, from definitions, the hypothesis implies that the lifted maps  $\Phi, \Phi': \widetilde M \to \widetilde B^n$ are so that $\Phi' \circ \Phi^{-1}$ is  the identity when restricted to the exceptional divisor $\pi^{-1}(0) = \bC P^{n-1}Ê$ and $\left(\Phi' \circ \Phi^{-1}\right)_*$ induces the identity map on each tangent space of a standard radial disc at the intersection with   $\bC P^{n-1} = \pi^{-1}(0)$. So, $Id_{\widetilde B^n}$ and $\Phi' \circ \Phi^{-1}$ are normalizing maps for $B^n$ that  satisfy the hypothesis of Theorem \ref{firstmaintheorem}. Hence $\Phi' \circ \Phi^{-1} = Id_{\widetilde B^n}$.\end{pf}
With the help of the previous lemma, we can prove the following.\par
\medskip
\begin{prop} \label{propparam} Fix a normalizing map $\Phi_o \in\cN(M, \tau, x_o)$ and a special frame  $(e^o_0, e^o_1, \dots, e^o_{n-1})$ at  $x_o\in M$,  associated with $\tau$. Also, for any $\Phi \in \cN(M, \tau, x_o)$  denote by $\phi: M \to B^n$ the corresponding projected map and let $(e_i^{\phi})$ be the  frame at $T_{x_o} M$ defined by $e_i^{\phi} = \left(\phi^{-1} \circ \phi_o\right)_*(e^o_i)$. \par
Then the frames of the form $(e_i^{\phi})$ are special frames, relatively to  $x_o$ and  $\tau$,  and  the correspondence $\Phi \longrightarrow (e_i^{\phi})$ is a  one to one map between $\cN(M, \tau, x_o)$ and the class of special frames at $x_o$.\end{prop}
\begin{pf} By definition of normalizing maps, the map $f = \phi \circ \phi_o^{-1}$ is a biholomorphism of $(M, J)$ into itself,  fixing $x_o$ and so that $\tau \circ f = \tau$.  This implies that $f_*|_{T_{x_o}M}$ maps the class of special frames at $x_o$ into itself.  Moreover, if $\Phi$ and $\Phi'$ are so that 
$(e_i^{\phi}) = (e_i^{\phi'})$, it follows from definitions that $\left.\left(\phi'\circ \phi^{-1}\right)_*\right|_{T_0 B^n} = Id$ and hence that $\Phi = \Phi'$ by the previous lemma. It remains to show that the correspondence $\Phi \longrightarrow (e_i^{\phi})$ is surjective. \par
Let $J^{(o)}$ be the complex structure on $B^n$ obtained by pushing forward the complex structure of  $M$ and denote by $(f^o_i)$ and $(f_i)$ the special frames of $(B^n, \tau_o, J^{(o)})$  obtained as images   of  $(e^o_i)$ and of another special frame  $(e_i)$, respectively.
 Let also $([f_i^o])$ and  $([f_i])$ the corresponding points in $\bC P^{n-1} \subset \widetilde B^n$. Observe  that the restriction  $\left.J^{(o)}\right|_{\bC P^{n-1} }$ coincides with the complex structure $J^{(o,0)}$ on $\widetilde B^n$ defined by the Bland and Duchamp invariant $\phi^{(0)}$ (see \S \ref{BD-invariant0}) and it is diffeomorphic to the standard complex structure of $\bC P^{n-1}$ (to see this, simply use the biholomorphism between $(\widetilde B^n, J^{(o,0)})$ and the blow up $\widetilde I$ of its indicatrix - see Proposition \ref{propositioncircular}).  Hence  $\operatorname{Aut}(\bC P^{n-1}, J^{(o)})$ is isomorphic to $\operatorname{Aut}(\bC P^{n-1}) = \operatorname{PGL}_{n}(\bC)$ and  there exists a unique $J^{(o)}$-biholomorphism $\Psi: \bC P^{n-1} \to \bC P^{n-1}$ mapping $[f_i]$ into $[f^o_i]$ for any $0 \leq i \leq n-1$. Let us  extend such a map to a diffeomorphisms $\Psi: \widetilde B^n \to \widetilde B^n$ in such a way that $\zeta \cdot \Psi(w) = \Psi(\zeta w)$ for any $\zeta \in \Delta$.
 We stress the fact that even if  $ \Psi|_{\bC P^{n-1}}$ is  by construction  a biholomorphism of $(\bC P^{n-1}, J^{(o)}|_{\bC P^{n-1}})$,  $\Psi$ is not  in general  a biholomorphism of $(\widetilde B^n, J^{(o)})$. \par
 Now, observe that $\Psi: \widetilde B^n \to B^n$ is a normalizing map which maps  the $[f_i]$'s into the $[f^o_i]$'s. Assume for the moment that $\Psi$ projects onto a map $\psi: B^n \to B^n$ which preserves the real manifold structure of $B^n$. Then,    $\psi_*(f_j) = e^{i \theta_j} f^{o}_j$ for some suitable complex numbers   $e^{i \theta_j} $ and, by a suitable adjustment of the definition of $\Psi$, we may always  assume that $e^{i \theta_j} = 1$ for any $0 \leq j \leq n-1$.   From this and its construction, it follows  that $\Phi = \Psi \circ \Phi_o$ is a normalizing map in  $\cN(M, \tau, x_o)$ so that $(e_i^{\phi}) = (e_i)$. This implies the  surjectivity of the map $\Phi \longrightarrow (e^{\phi}_i)$.\par
 So, in order to conclude, we only need to show that $\Psi$ projects onto a map $\psi: B^n \to B^n$ which preserves the real manifold structure of $B^n$. This is equivalent to check that there is a chart on $B^n$ on a neighborhood of  the origin, which belongs to the real manifold structure of $(B^n, J^{(o)})$ and in which $\psi$ is smooth. This fact can be done  using the circular representation. In fact, we may use it  in order to identify the differentiable manifolds $(B^n, J^{(o)})$ and $(\widetilde B^n, J^{(o)})$ with the indicatrix $I$ and its  blow up $\widetilde I \subset \widetilde \bC^n$,  respectively, both endowed with a suitable  complex structure, say  $\hat J^{(o)}$. Since the circular representation is a diffeomorphism between a domain in $\bC^n$ and the manifold of circular type,  we have that  the real manifold structure  on $I$ determined by projection from the manifold structure of $\widetilde I$ is the same of  the standard manifold structure  of $I$, considered as open subset of $\bC^n$. Writing the explicit expression   $\Psi$ as a map $\Psi: \widetilde I \to \widetilde I$ it can be checked directly that the projected map  $\psi: I \to I$   is smooth in any system of standard coordinates of $\bC^n$. \end{pf}
The  main result of the section can be now  immediately inferred. \par
\medskip
\begin{theo} \label{parameterizationnormalizingmaps} Let $(D, J, \tau)$ be a stable domain of circular type and $\cN(D)$ the class of all normalizing maps, which induce on $B^n$ the same  real manifold structure. Fix also a frame  $(e^o_i) \subset T_{x} D$ belonging to the pseudo-bundle of special frames $\pi: P\to D$ and a normalizing map $\Phi_o \in \cN(D)$. Then, the correspondence $\Phi \longrightarrow (e_i^{\phi})$ of the previous proposition gives a continuous one to one map between $\cN(M)$ and the pseudo-bundle of special frames $P$.\end{theo}
We remark that if $P$ is identified with a subset of $\operatorname{U}_n(B^n) \simeq \operatorname{Aut}(B^n)$, the previous result can be stated saying that: {\it $\cN(D)$ is parameterized by a special subset of  $  \operatorname{Aut}(B^n)$\/}. In particular, we have that if  the set of centers $\mathfrak C(D)$ is a singleton,  then $\cN(D)$ is parameterized by  $\operatorname{U}_n =  \operatorname{Aut}(B^n)_0$, while if $\mathfrak C(D) = D$,  $\cN(D)$ is parameterized by  $ \operatorname{Aut}(B^n)$. It is interesting to observe the analogy (and the difference) between this class of  normalizing maps and  the class of Chern-Moser normalizing maps for Levi non-degenerate hypersurfaces of $\bC^n$. \par\bigskip
\subsubsection{The automorphism group of a manifold of circular type}\hfill\par
\medskip
Let $M$ be a circular domain of circular type. By means of a normalizing map, there is no loss of  generality if we assume that  $M$ is  $(B^n, J, \tau_o)$, where $J$ is a complex structure of Lempert type, L-isotopic to the standard one. In this case, any  automorphism  $\Phi \in \operatorname{Aut}(M) = \operatorname{Aut}(B^n, J)$ is also a normalizing map and, if  we set $\Phi_o = Id_{B^n}$, Theorem \ref{parameterizationnormalizingmaps} implies that the action of $G = \operatorname{Aut}(M)$ on special  frames determines a one to one map between $G$ and the points of any  orbits  $G \cdot (e^o_i)$.  In particular, if $M$ is a stable domain of circular type and $G$ is a Lie group, then  $G$ is diffeomorphic to any of its orbits $G \cdot u$ in  the pseudo-bundle $P$. In such a case,   if   $\mathfrak C(M) \subset M$  is contained in a real submanifold of dimension $a$ ($\leq 2n = \dim_\bR M$), we get that 
$$\dim G \leq a + n^2 \leq 2n + n^2\ .$$\par
By the classical results of H. Cartan (see e.g. \cite{Ak}), the property that $G$ is a Lie group is granted whenever $M$ is  biholomorphic to a bounded complex domain in $\bC^n$. The previous remark gives therefore  a refinement of the well known upper bound $2n + n^2$  for $\dim G$. \par
We  expect that   $G = \operatorname{Aut}(M)$ is a Lie group for {\it any} stable domain (not necessarily embedded in $\bC^n$) and hence that the above estimate is true in full generality: We plan to attack such  a problem in a future paper.\par
\par
\bigskip
\subsection{Characterizations of circular domains and of the unit ball}\hfill\par
\bigskip
\begin{definition}{\rm Let $(D, J, \tau)$ be a stable domain of circular type. 
An element $g \in G = \operatorname{Aut}(D,J)$ is called {\it rotational at the center $x_o$\/}  if 
$g_*|_{x_o} = \lambda Id_{T_{x_o} D}$ for some $\lambda \in \bC^n$ so that $\lambda^n \neq 1$ for any $n \in \bZ$. We say that {\it $M$ is rotational at $x_o$\/} if $G$ contains a rotational element at $x_o$.}\end{definition}
The following two theorems are extensions to the case of stable domains of a couple of results in \cite{Pt}.\par
\begin{theo} A stable domain of circular type $(D, J, \tau)$ is biholomorphic to a circular domain in $\bC^n$ if and only if it is rotational at some point $x_o \in D$. 
\end{theo}
\begin{pf}  There is no loss of generality if we assume that $D$ is in normal form, i.e. that $(D, J, \tau) = (B^n, J, \tau_o)$, and that $x_o = 0$. If $g$ is rotational, the lifted map $\widetilde g: \widetilde B^n \to \widetilde B^n$ is so that $\widetilde g|_{\bC P^{n-1}} = Id_{\bC P^{n-1}}$. Moreover, by Lemma \ref{uniqueness},  $\tau_o \circ \widetilde g = \tau_o$ and hence  $\widetilde g$ induces on any standard radial disc a biholomorphism  that fixes the origin and with derivative at $0$ equal to $\lambda$. This implies that $\lambda = e^{i \theta}$ for some $\theta \in \bR$ and,  in a system of coordinates as described before (2.5), $\widetilde g$ is of the form $\widetilde g(v, \zeta) = (v,  e^{i \theta} \zeta)$. Moreover, since $\widetilde g$ is a $J$-biholomorphism and induces on $\bC P^{n-1}$ the identity map, all Bland and Duchamp invariants do not change under the action of $\widetilde g$. In particular any invariant $\phi^{(k)}$, $k \geq 1$, is so that 
$\phi^{(k)} -  e^{i k \theta}\phi^{(k)} = 0$.
Since $\lambda^k = e^{i k \theta} \neq 1$ for all $k$, we get  that $\phi^{(k)} = 0$ for any $k \geq 1$ and the conclusion follows from Proposition \ref{propositioncircular} b). 
\end{pf}
\begin{theo} \label{twocenters} A stable domain $(D, J, \tau)$ is biholomorphic to the unit ball $B^n$ if and only if it  is rotational at two distinct centers $x_o, x_o' \in D$. 
\end{theo}
\begin{pf} The necessity is trivial. For the sufficiency notice that,  by the previous theorem, $D$ is biholomorphic to a circular domain  $I \subset \bC^n$, which is rotational w.r.t. to two distinct points. The conclusion follows automatically from \cite{Pt}, Thm. 9.4.
\end{pf}
In the following statement, for a given manifold of circular type $(M, J, \tau)$ with $\tau: M \to [0, r^2)$ and center $x_o$, for any $c  \in (0,r^2)$, we denote by $D_{< c}$ the domain contained in $M$ defined by $M = \{\ x \in M\ , \ \tau(x) < c\ \}$. We remark that the indicatrix at $x_o$ of $(M, J, \tau)$ coincides (up to rescaling) with the indicatrix at that point of $(D_{<c}, J, \tau)$. From this and   Proposition \ref{circularrepresentation}, 
it can be directly checked that  $(D_{<c}, J, \tau)$  is always a  stable domain if $D_{<c} \subsetneq M$. \par
\begin{theo} \label{lasttheorem} Let $(M, J, \tau)$ be a manifold of circular type of dimension $n$.  Then $M$ is biholomorphic to a circular domain in $\bC^n$ if and only if for some given parabolic exhaustive function $\tau$ there exist  
  two domains  $D_{<c}$, $D_{<{c'}}$,  $0 < c < c' < r^2$, that are  biholomorphic one to the other by means of a map fixing the center of $\tau$.
\end{theo}
\begin{pf} The necessity is direct: If $M \subset \bC^n$ is a circular domain with Minkowski function $\mu$ and we set $\tau = \mu^2$, then the map $f(z) = \frac{c'}{c} z$ determines a biholomorphism between $D_{<c}$ and $D_{c'}$ for any two $c,c'$. \par
To prove the sufficiency, we may  assume that $M$ is in normal form, i.e. $M = (B^n, J, \tau_o)$ and that 
the parabolic exhaustive function which defines the two biholomorphic stable domains $D_{<c}$ and $D_{<c'}$ is the function $\tau = \tau_o$. Notice also that the restriction on $D_{<c'}$ of the map
$$\psi: \widetilde \bC^n \to \widetilde \bC^n \, \qquad \psi([z], z) = ([z], \frac{1}{c'} z)$$
is a normalizing map for $D_{<c'}$ (i.e. maps $(D_{<c'}, J, \tau_o)$ into $(B^n, \tilde J, \tau_o)$) and maps $D_{<c}$ into $D_k$, $k = c/c'$.  \par
It follows directly from  definitions that  the indicatrix $I$ at $x_o = 0$ of  $(B^n, \tilde J, \tau_o)$ is the same of the  $(D_k, \tilde J, \tau_o)$ (up to rescaling) and hence all special frames of the first domain at $0$ coincide with the special frames of the second domain up to multiplication by $k$. On the other hand, by hypothesis and Lemma 2.2, we have a biholomorphism of domains of circular type $f:   (D_{<k}, \tilde J, \tau_o)\to (B^n, \tilde J, \tau_o) $ so that $f(0) = 0$. The differential $f_*|_0:  T_0 D_{<k} = T_0 B^n \to T_0 B^n  $ is a $\bC$-linear map mapping  $k \cdot I$ into $I$ and hence mapping a fixed special frame $(e_i)$  into another special frame $(e'_i)$  rescaled by the factor $1/k$. Let us denote by $\hat f: B^n \to B^n$ the unique normalizing map of $(B^n, \tilde J, \tau_o)$ which transforms $(e_i)$ into $(e'_i)$ (see Theorem \ref{parameterizationnormalizingmaps}).  The lift at the blow up level of  $f^{-1} \circ \hat f$  is a diffeomorphism  between $\widetilde B^n$ and $\widetilde D_{<k} \subset \widetilde B^n$ which induces the identity map on $\bC P^{n-1}$ and so that, when restricted to any radial disc, it is a holomorphic map (w.r.t. the standard complex structure) fixing the center, mapping the unit disc  into the disc of radius $k$ and with derivative equal to $k$ at the origin. Therefore $f^{-1} \circ \hat f$  is of the form 
$$ (f^{-1} \circ \hat f)([z], z)    =\psi_k([z], z)\ ,\ \ \text{where}\ \ \psi_k([z], z) \= ([z], kz)\ .$$
The same argument can be repeated for any iterated map $f^n ={f \circ \dots \circ f}$ and we obtain that, for any $n \in \mathbb N$, the map 
$$\hat f^{(n)} = f^n \circ \psi_k^n$$
is a normalizing map of $(B^n, \tilde J, \tau_o)$, fixing the origin.  Since the normalizing maps fixing the origin are continuously parameterized by the compact set of special frames at the origin ($\simeq\operatorname{U}_n$), we may consider a subsequence $n_j$ so that the sequence of normalizing maps $\hat f^{(n_j)}$ converges  uniformly on compacta to a normalizing map $\hat f^{(\infty)}$. In particular,   the sequence of 
complex structures $\tilde J^{(n_j)} \= \hat f^{(n_j)}{}^*(\tilde J)$ converges to a complex structure $\tilde J^{(\infty)} =  \hat f^{(\infty)}{}^*(\tilde J)$. On the other hand, since $f$ is a $(\tilde J, \tilde J)$-biholomorphism, we have that 
$$\tilde J^{(n_j)}  = \hat f^{(n_j)}{}^*(\tilde J) = \psi_k^{n_j}{}^*(\tilde J)\ .$$
A direct computation shows that, for any point $([z], z) \in \widetilde B^n$,  the deformation tensor $\phi^{(n_j)}|_{([z],z)}$ of   $\psi_k^{n_j}{}^*(\tilde J)$ coincides with the deformation tensor $\phi$ of $\tilde J$, but evaluated at the point $([z], k^{n_j} z)$. It follows that 
$$\phi^{(\infty)}|_{([z], z)} = \lim_{n_j \to \infty} \phi^{(n_j)}|_{([z],z)}  = \lim_{n_j \to \infty} \phi|_{([z], k^{n_j} z)} =  \phi|_{([z], 0)} \ ,$$
i.e. the deformation tensor $\phi^{(\infty)}$ of $\widetilde J^{(\infty)}$ is independent on the coordinate of the radial discs.  By Proposition \ref{propositioncircular}, this means that the domains
$$(B^n, \widetilde J^{(\infty)}, \tau_o) \simeq (B^n, \widetilde J, \tau_o) \simeq (D_{<c'},  J, \tau_o) $$
are all circular. In particular, the restriction to $D_{<c'}$ of the deformation tensor of the manifold from  which started, i.e. of $(B^n, J, \tau_o)$,  is independent on the coordinates of the radial discs.  By Proposition \ref{invariants} (iii) (i.e. analyticity along the radial discs),  we get 
independence on the coordinates of the radial discs  over the entire  $(B^n, J, \tau_o)$. Using once again   Proposition \ref{propositioncircular}, we get the claim, i.e. $(B^n, J, \tau_o)$ is circular.
\end{pf}
\medskip
\begin{rem}{\rm Notice that  the proof of previous theorem implies that,  if  we assume that $(M, J, \tau)$ is  a stable domain (in particular if it is  a strictly convex domain in $\bC^n$), the claim is true also when $c' = r^2$ and  $D_{<c'} = M$. }
\end{rem}
By the previous theorem and remark and by Theorem \ref{twocenters}, we have the following direct corollary.\par
\begin{cor} \label{lastcorollary}ÊA stable domain $D$ is biholomorphic to the unit ball if and only it has at least   two centers $x_o\neq x'_o$ and  it is biholomorphic, by means of two maps fixing $x_o$ and $x_o'$, respectively,  to two   proper subdomains $D_{<c}=\{ \ \tau < c\ \}$, $D'_{<c'}=\{ \ \tau' < c'\ \}$, with   $\tau$, $\tau'$ parabolic exhaustive functions  associated with $x_o$ and $x'_o$, respectively.
\end{cor}

\bigskip
\bigskip
\font\smallsmc = cmcsc8
\font\smalltt = cmtt8
\font\smallit = cmti8
\hbox{\parindent=0pt\parskip=0pt
\vbox{\baselineskip 9.5 pt \hsize=3.1truein
\obeylines
{\smallsmc
Giorgio Patrizio
Dip. Matematica ``U. Dini''
Universit\`a di Firenze
Viale Morgani 67/a
I-50134 Firenze
ITALY
}\medskip
{\smallit E-mail}\/: {\smalltt patrizio@math.unifi.it
}
}
\hskip 0.0truecm
\vbox{\baselineskip 9.5 pt \hsize=3.7truein
\obeylines
{\smallsmc
Andrea Spiro
Dip. Matematica e Informatica
Universit\`a di Camerino
Via Madonna delle Carceri
I-62032 Camerino (Macerata)
ITALY
}\medskip
{\smallit E-mail}\/: {\smalltt andrea.spiro@unicam.it}
}
}

\end{document}